\documentclass[11pt,a4paper,leqno,twoside]{amsart}
\usepackage{longtable}

\newcommand{\qz}{{\mathbb Q}}

\newcommand{\pz}{{\mathbb P}}
\newcommand{\fz}{{\mathbb F}}

\begin{document}

\title{Modularity experiments on $S_4$-symmetric double octics}
\author{Christian Meyer}
\email{cm.math@gmx.de}
\dedicatory{To Duco van Straten, on the occasion of his 60th birthday}
\date{October 10, 2018}

\begin{abstract}
We will invest quite some computer power to find double octic threefolds that are connected to weight four modular forms.
\end{abstract}

\maketitle

\section{Introduction}

More than twenty\footnote{Previous versions of this paper had 'ten' in this place. That's life.} years ago, under the supervision of Duco van Straten in Mainz, I wrote my doctoral thesis, which finally made it into book form [0]. Shortly afterwards, I left academia for rocket science\footnote{That's how they describe what mathematicians are doing at banks.}\textsuperscript{,}\footnote{I don't like footnotes, but I will make exceptions for Sir Terry Pratchett, David Foster Wallace, and this paper.}. However, I could not stand my desktop computers running idle much of the time, so I took up modular Calabi-Yau threefolds, the topic of my thesis, again. It turned out rather convenient: with just the effort of a few days of programming, I could lean back (or perform my day job) and let the machines produce interesting results, providing the illusion of still doing mathematics.

In the following I will describe some of the results. The paper is not meant for 'serious' publication\footnote{The arXiv IS serious, of course, but not 'serious', since you do not have to pay for access.}, which grants me some freedom for un-academical rambling. Nevertheless, many interesting examples might be hidden in the mess, and current or future students of mathematics are invited to take a closer look. I will also use the opportunity to provide updated references for my book.

\section{$S_4$-symmetric double octics}

A double octic is a double covering of $\pz^3$ branched along an octic surface $D$. It can be represented as a hypersurface in the weighted projective space $\pz^4(4,1,1,1,1)$, with equation
\[
u^2 = f(x:y:z:t)
\]
where $f$ is a homogeneous polynomial of degree eight. Singularities of $X$ correspond to singularities of $D$ and can be resolved by blow-up. Depending on the types of the singularities of $D$, the resulting threefold may be Calabi-Yau or not, but arithmetically this does not make too much of a difference.

In order to find interesting examples (for instance, with many singularities of a certain kind, or interesting arithmetical properties) it is a useful strategy to have a large group acting on $D$. We will be concerned with the case of octics that are invariant under action of the symmetric group $S_4$ by permutation of coordinates.

The space of these $S_4$-symmetric octics is generated by the degree eight products of the elementary-symmetric polynomials on four variables,
\begin{align*}
e_4 := e_4(x,y,z,t) & = xyzt,\\
e_3 := e_3(x,y,z,t) & = xyz + xyt + xzt + yzt,\\
e_2 := e_2(x,y,z,t) & = xy + xz + xt + yz + yt + zt,\\
e_1 := e_1(x,y,z,t) & = x + y + z + t.
\end{align*}
There are 15 of them, and in order to save on indices, I assigned them to a 15-letter heterogram, BRUCHWEGSTADION
\footnote{The stadium in which the Mainz 05 football team used to play until 2011. Its successor started as 'Coface Arena', was renamed to 'Opel Arena', and now runs under 'Mewa Arena'. I propose to name all stadiums 'Mammon Arena [insert club]'.}\textsuperscript{,}\footnote{Another nice one would have been UNPROBLEMATISCH. And then WHISKYPRODUZENT, which would also have worked, with the Z replaced by C, in Dutch.}.
\[
\begin{array}{lll}
B : e_4^2         & \qquad W : e_3^2 e_2     & \qquad A : e_2^4\\
R : e_4 e_3 e_1   & \qquad E : e_3^2 e_1^2   & \qquad D : e_2^3 e_1^2\\
U : e_4 e_2^2     & \qquad G : e_3 e_2^2 e_1 & \qquad I : e_2^2 e_1^4\\
C : e_4 e_2 e_1^2 & \qquad S : e_3 e_2 e_1^3 & \qquad O : e_2 e_1^6\\
H : e_4 e_1^4     & \qquad T : e_3 e_1^5     & \qquad N : e_1^8
\end{array}
\]
If the coefficients of an octic are defined over $\qz$, both the surface and the corresponding double octic are also defined over $\fz_p$, the finite field with $p$ elements for a prime number $p$.

\section{Counting points}
\label{sec_counting}

Modularity refers to the phenomenon that the number of points on Calabi-Yau manifolds over finite fields seems to be related to the coefficients of certain modular forms. Specifically, in dimension 3, one might find relations of the form
\[
a_p \equiv 1 - \# X_p \quad \mod p
\]
for almost all primes $p$, where $a_p$ are the coefficients of a newform of weight 4 for $\Gamma_0(N)$, and $\# X_p$ are the numbers of points on (a singular model of) the threefold $X$ over $\fz_p$. The relation should hold for all primes except possibly the divisors of the level $N$ (the \emph{bad} primes).

How to find these rare birds? In this paper we will be using brute force. In order to count points on an $S4$-symmetric double octic
\[
u^2 = f(x:y:z:t)
\]
we may plug all points $(x:y:z:t)\in\pz^3(\fz_p)$ into $f$, and count one if $f(x:y:z:t)=0$, and two if $f(x:y:z:t)$ is a square modulo $p$. We can then substantially speed up the process (1) by pre-computing the values of the BRUCHWEGSTADION polynomials, and (2) by aggregating points that lead to the same set of values. The aggregates may correspond to $S_4$ orbits but may also be larger, particularly for small $p$. For the first 25 primes (up to 97), the number of aggregates amounts to $759,931$, only $\approx 16\%$ of $4,763,331$, the sum over the numbers of points in $\pz^3(\fz_p)$. Using these short cuts it is possible to process many examples in relatively short time. More on this in Section \ref{sec_usefulness}.

\section{Correspondences}

Sometimes the same modular form can be linked to different threefolds. In this case, by the Tate conjecture, there should be a correspondence between the threefolds, i.e., an algebraic cycle on the product of the two varieties.\footnote{The two threefolds will then be called \emph{relatives}.} Finding correspondences (e.g., finite maps between two threefolds) is highly non-trivial since one needs an \emph{idea} first. On the other hand, once a correspondence machine has been found, new modularity examples (with maybe very different geometry) can be generated. We will explain some of these machines for our case of $S_4$-symmetric double octics\footnote{As Dave Barry might say, \emph{The Double Octics} would be an excellent name for a rock band. And \emph{The Correspondence Machine} might not be too bad for the prog variety.}, and save some space in the tables by removing much that can be explained by correspondence.\footnote{In the end, only one example per modular form should remain, of course. Unfortunately, I ran out of \emph{ideas}.}

There are at least four types of correspondences that work for some or all $S_4$-symmetric double octics:

\subsection{Coordinate change}
\label{sec_coord}

The map
\begin{align*}
(x:y:z:t) \mapsto (S_1 + \lambda x:S_1 +\lambda y:S_1 + \lambda z:S_1 + \lambda t)
\end{align*}
with $(1:\lambda)\in\pz^1$ treats the elementary-symmetric polynomials as follows:
\begin{align*}
e_1 & \mapsto (4+\lambda) e_1\\
e_2 & \mapsto \lambda^2 e_2 + 3(2+\lambda) e_1^2\\
e_3 & \mapsto \lambda^3 e_3 + 2\lambda^2 e_2 e_1 + (4+3\lambda) e_1^3\\
e_4 & \mapsto \lambda^4 e_4 + \lambda^3 e_3 e_1 + \lambda^2 e_2 e_1^2 + (1+\lambda) e_1^4
\end{align*}
For $\lambda\neq -4$ this map induces a correspondence but for almost all choices of $\lambda$, the resulting coefficients of the BRUCHWEGSTADION polynomials will be too large to already have been found by the brute force search. Exceptions will be marked in the table of results.

\subsection{The Segre map}
\label{sec_segre}

The Segre map
\[
(x:y:z:t) \mapsto (x^2:y^2:z^2:t^2)
\]
treats the elementary-symmetric polynomials as follows:
\begin{align*}
e_1 & \mapsto e_1^2 - 2e_2\\
e_2 & \mapsto e_2^2 + 2e_4 - e_3 e_1\\
e_3 & \mapsto e_3^2 - 2e_4 e_2\\
e_4 & \mapsto e_4^2
\end{align*}
In particular, for those octics exhibiting an $e_4$ factor (i.e., the $BRUCH$ ones) that factor will be squared, thus eliminated in a double covering, and a correspondence will be induced as follows:
\begin{align*}
B & \mapsto B\\
R & \mapsto 4U - 2C -2W + E\\
U & \mapsto 4B - 8R + 4U + 4E - 4G + A\\
C & \mapsto 8U - 8C + 2H - 8G + 8S - 2T + 4A - 4D + I\\
H & \mapsto 16A - 32D + 24I - 8O + N
\end{align*}
Not the most elegant of relationships perhaps, but it shows that for each of the $BRUCHWEGSTADION$ polynomials there are modular examples with non-zero coefficients for this polynomial -- and that they may be hard to find by our brute force approach, due to restricting relations between the polynomials that are driving up the coefficients. 

\subsection{Inversion}
\label{sec_BRUCWEGA}

The map
\[
(x:y:z:t) \mapsto \left(\frac{e_4}{x}:\frac{e_4}{y}:\frac{e_4}{z}:\frac{e_4}{t}\right) = (yzt:xzt:xyt:xyz)
\]
treats the elementary-symmetric polynomials as follows:
\begin{align*}
e_1 & \mapsto e_3\\
e_2 & \mapsto e_4 e_2\\
e_3 & \mapsto e_4^2 e_1\\
e_4 & \mapsto e_4^3
\end{align*}
For some of the $BRUCHWEGSTADION$ polynomials, a fourth power of $e_4$ will appear that can be eliminated from the equation for the double covering. The resulting polynomial has degree eight again. The polynomials in question are the $BRUCWEGA$ ones, and the only thing to happen under the map is $C$ and $W$ changing places:
\begin{align*}
B & \mapsto B\\
R & \mapsto R\\
U & \mapsto U\\
C & \mapsto W\\
W & \mapsto C\\
E & \mapsto E\\
G & \mapsto G\\
A & \mapsto A
\end{align*}
In contrast to the Segre-type correspondences, these are visible to the naked eye.

\subsection{Sign change}
\label{sec_sign}

The map
\[
(x:y:z:t) \mapsto (x:y:z:-t)
\]
maps some $S_4$-symmetric polynomials to $S_4$-symmetric polynomials again. In particular, with
\[
P := 8 e_3 - 4 e_2 e_1 + e_1^3 = (x-y-z+t)(x-y+z-t)(x+y-z-t)
\]
we find
\[
e_1 P \mapsto e_1 P - 16 e_4.
\]
Considering additionally the obvious invariant $e_1^2-2 e_2 = x^2+y^2+z^2+t^2$ we find
\begin{align*}
B = e_4 & \mapsto -B\\
8 R - 4 C + H = e_1 P & \mapsto -16 B + 8 R - 4 C + H\\
4 U - 4 C + H = (e_1^2-2 e_2)^2 & \mapsto 4 U - 4 C + H
\end{align*}
and, after some scaling,
\begin{align*}
\alpha\cdot B + 2\beta \cdot R + (4\gamma - \beta)\cdot U - 4\gamma\cdot C + \gamma\cdot H\\
\mapsto
-(4\beta + \alpha) \cdot B + 2\beta\cdot R + (4\gamma-\beta)\cdot U - 4\gamma\cdot C + \gamma\cdot H
\end{align*}
There will be additional correspondences outside BRUCH, of course, but they are less relevant for the subset of examples to be considered.

\section{Usefulness}
\label{sec_usefulness}

Armed with the mathematical apparatus described above, I had written a C++ program for counting points on $S_4$-symmetric double octics. Using OpenMP for parallelization, the program is able to check more than one million examples in less than two minutes on the CPU at my disposal, an AMD Ryzen 9 3900 with twelve cores and hyperthreading (which replaced an Intel Xeon E3-1230 V2 with four cores and hyperthreading, in 2020). Curiously, this is still a '3 Gigahertz machine', but around 30 times faster\footnote{Or probably more; this is difficult for me to estimate since the CPU databases available on the internet are less informative about such veterans, and I do not have a Pentium IV machine at hand, so I cannot run my own benchmarks.}\textsuperscript{,}\footnote{And vastly more efficient!} than the Pentium IV described as such in [0], Section 1.8.

The table of newforms of weight 4 for $\Gamma_0(N)$ was an extension of the one from [0]: 4033 newforms altogether, complete up to level $N=2500$, with $N=2700$ added. The table contains the coefficients for the first 25 primes (up to 97). A double octic will be shown in the table if the relation from Section \ref{sec_counting} holds for at least 21 primes. In order to always find the twist of minimal level (cf. [0], Section 1.7.3), the program also modifies the signs of the coefficients.

I did not track the CPU time spent on the experiment but a lower bound of five years should be safe to assume. With around 50W of extra power consumption (compared to the computer running idle) this amounts to an impressive $2.2$MWh.\footnote{Of course this is nothing in comparison with the insane amounts of energy spent on the bitcoin delusion and other blockchain endeavours, and more recently the stochastic parrots called artificial intelligence.}

At least we should reassure ourselves that I have not been fooled by randomness\footnote{This expression has been coined by N.N. Taleb, cf. {\tt www.fooledbyrandomness.com}. If you haven't heard about him, put away this paper for later and have a look. Unfortunately, he went a little crazy over Covid.}. Wouldn't we expect to find (at least some) spurious coincidences if we have been searching for long enough?

Let's do a back-of-the-envelope computation. Given the 25 coefficients (modulo the respective primes) of a newform from the table, what is the probability of hitting at least 21 of them by drawing randomly from $\fz_2\times \fz_3 \times \ldots\fz_{97}$? Consider the set
\[
\tilde{P} = \{1, 2, 4, 6, 10, \ldots, 96\}
\]
of the first 25 primes reduced by one. Hitting at least 21 coefficients means missing at most four, and there are $e_n(\tilde{P}^{25})$ possibilities where here $e_n$ is the elementary-symmetric polynomial of degree $n$ in 25 variables. Evaluation is not too difficult (e.g., by exploiting the relations between the $e_n$ and the power sums in 25 variables\footnote{At this point, you might switch to a larger envelope, or to a spreadsheet.}), and we find that in total there are
\[
1 + e_1(\tilde{P}^{25}) + e_2(\tilde{P}^{25}) + e_3(\tilde{P}^{25}) + e_4(\tilde{P}^{25}) \approx 3.3\cdot 10^{11}
\]
possibilities for hitting at least 21 coefficients. On the other hand, the number of points in $\fz_2\times \fz_3 \times \ldots\fz_{97}$ is around $2.3\cdot 10^{36}$. Altogether, assuming that the coefficients of the newforms can be interpreted as the results of random draws\footnote{Which they can not, of course.}, and that the numbers of points on the double octics can be interpreted as the results of random draws\footnote{Which they can not, of course.}, the expected number of chance matchings is of the order of magnitude
\[
\#\text{octics tried}\cdot \#\text{newforms}\cdot \frac{3.3\cdot 10^{11}}{2.3\cdot 10^{36}}
\approx \#\text{octics tried}\cdot 4.9\cdot 10^{-23}.
\]
Therefore, as long as the number of octics tried stays some orders of magnitude below $10^{23}$ (as it does, cf. the following section) we do not have to worry too much about chance occurrences.

\section{Results}
\label{sec_tables}

Now it is time for actual results. I concentrated the search on small absolute values of the coefficients. Let $b,\ldots,n$ denote the coefficients of the polynomials $B,\ldots,N$,
\[
\Psi := |b|+|r|+|u|+|c|+|h|+|w|+|e|+|g|+|s|+|t|+|a|+|d|+|i|+|o|+|n|
\]
and $\Phi$ the number of non-zero coefficients among $b,\ldots,n$. The search ran at least\footnote{The careful reader will detect some additional examples in the table.} over the following sets of parameters (leading to around $10^{12}$ double octics to be examined):
\begin{itemize}
\item $\Phi\leq 2$, $\Psi\leq 2000$
\item $\Phi = 3$, $\Psi\leq 500$
\item $\Phi = 4$, $\Psi\leq 180$
\item $\Phi = 5$, $\Psi\leq 60$
\item $\Phi = 6$, $\Psi\leq 35$
\item $\Phi = 7$, $\Psi\leq 25$
\item $\Phi\geq 8$, $\Psi\leq 18$
\item BRUCH only (cf. Section \ref{sec_segre}), $\Psi\leq 270$
\item BRUCWEGA only (cf. Section \ref{sec_BRUCWEGA}), $\Psi\leq 50$
\end{itemize}
This way I found many examples -- and surely missed many interesting ones. The correspondences presented in Sections \ref{sec_coord} and \ref{sec_segre}, at least, will usually lead to examples with both larger $\Psi$ and larger $\Phi$. On the other hand, known examples with larger symmetry group (e.g., $S_5$-symmetric octics or Heisenberg-invariant octics, cf. [0], Sections 4.10 and 4.11) will typically exhibit $\Psi$ and $\Phi$ too large to be rediscovered.

In the following table we will collect the coefficients $b,\ldots,o$ ($n$ will be omitted since it always equals zero\footnote{In the examples found by the search, that is. By applying a correspondence machine, examples with nonzero $n$ can be generated.}), the associated newform using the notation from [0], and comments about applicability of some correspondence machine(s). Correspondences will be marked as follows:
\begin{itemize}
\item $S$: a relative can be found by applying the Segre construction from Section \ref{sec_segre}. This relative will not be listed in the table.
\item $I$: a relative can be found by applying the inversion construction from Section \ref{sec_BRUCWEGA}. This relative will not be listed in the table.
\item $L_1, L_2, L_3$: two threefolds marked by $L_i$ are related by a coordinate change as in Section \ref{sec_coord}. The parameter $\lambda$ amounts to -2 for $L_1$, and to -1 for $L_2$ and $L_3$.
\item $C_1,\ldots,C_7$: two threefolds marked by $C_i$ are related by sign change as in Section \ref{sec_sign}.
\item $X_1, X_2$: two threefolds marked by $X_i$ are related by the coordinate change $x\mapsto x+y$, $y\mapsto x-y$, $z\mapsto z+t$, $t\mapsto z-t$, cf. [0], Section 4.11.3.
\end{itemize}
References to [0] and other sources, as well as specific comments, will be relegated to footnotes.

The modular forms detected involve the following bad primes (i.e., level divisors): 2, 3, 5, 7, 11, 13, 17, 19, 23, 29, 31, 37, 43, 47, 53, 67, 73, 101, 103, 113. Some gaps in the list of modular forms connected to threefolds are closed, e.g., for the newform 23/1. However, I do not consider my brute force approach to be fruitful beyond this point. For \emph{really} interesting examples of modular threefolds, including the one that first filled the gap at 7/1 (there's another one in the table below), I recommend [114] and [138].

\setlength{\tabcolsep}{4.5pt}
\begin{center}
\begin{scriptsize}
\begin{longtable}{rrrrrrrrrrrrrrll}
\hline
$b$ & $r$ & $u$ & $c$ & $h$ & $w$ & $e$ & $g$ & $s$ & $t$ & $a$ & $d$ & $i$ & $o$ & form & corr.\\
\hline
\endhead
\hline
\endfoot
 &  &  &  &  &  & 64 & -64 & 16 & -1 &  &  &  &  & 5/1 & \\
 &  &  & 4 & -2 & -4 & -2 & 4 &  &  & -1 &  &  &  & 5/1 & \\
 &  & 4 & -4 &  &  & 1 & -2 &  &  & 1 &  &  &  & 5/1 & $I$\\
 & 8 & -4 &  &  &  &  & -2 &  &  & 1 &  &  &  & 5/1 & \\
 & 8 & 2 & 1 &  & -3 &  &  &  &  &  &  &  &  & 5/1 & $I$\\
16 & 4 & -8 & 6 & & -2 & & -3 & & & 1 & & & & 5/1 & $I$\\
16 & 4 & -4 & 1 &  & 1 &  &  &  &  &  &  &  &  & 5/1 & \\
256 & -32 &  &  & 1 &  &  &  &  &  &  &  &  &  & 5/1 & $S$\\
\hline
 &  &  &  &  &  &  & 16 & 2 &  &  & -24 & -3 &  & 6/1 & \\
 &  &  &  &  &  &  & 48 & -18 &  &  & -8 & 3 &  & 6/1 & \\
 &  &  &  &  &  &  & 64 & -54 &  &  & -32 & 25 &  & 6/1 \\
 &  &  &  &  &  &  & 128 & -68 & 9 &  &  &  &  & 6/1 & \\
 &  &  &  &  &  & 18 & -17 &  &  & 4 &  &  &  & 6/1 & \\
 &  &  &  &  & 8 &  & -16 & 2 &  &  & 6 & -3 &  & 6/1 & \\
 &  &  &  &  & 8 &  &  & 2 &  &  & -2 & -1 &  & 6/1 & \\
 &  &  &  &  & 8 &  & 24 & 2 &  &  & 10 & 5 &  & 6/1 & \\
 &  &  &  &  & 48 &  & -8 &  &  & -8 & 3 &  &  & 6/1 & \\
 &  &  & 1 &  & 1 &  & -3 &  &  &  &  &  &  & 6/1 & \\
 &  &  & 1 &  & 1 &  &  &  &  &  &  &  &  & 6/1 & \\
 &  &  & 2 &  &  &  & -1 &  &  &  &  &  &  & 6/1 & $I$\\
 &  &  & 2 &  & 2 &  & 2 &  &  &  & -4 & 1 &  & 6/1 & \\
 &  &  & 4 &  & 4 &  & -8 &  &  &  & 4 & -1 &  & 6/1 & \\
 &  &  & 8 & -3 & -16 & 6 &  &  &  &  &  &  &  & 6/1 & \\
 &  &  & 12 &  & -12 &  & 8 &  &  &  & -4 & 1 &  & 6/1 & \\
 &  &  & 12 &  &  & -9 & 2 &  &  &  &  &  &  & 6/1 & $I$\\
 &  &  & 48 &  &  & -9 & -4 &  &  & -4 &  &  &  & 6/1 & $I$\\
 &  &  & 54 &  &  &  & -27 &  &  & 8 &  &  &  & 6/1 & $I$\\
 &  & 8 & -9 &  & -9 & -12 & 7 &  &  & -2 &  &  &  & 6/1 & $I$\\
 &  & 8 & -9 &  & -9 &  & 7 &  &  & -2 &  &  &  & 6/1 & \\
 &  & 8 & -3 &  & -3 & -9 & 9 &  &  & -2 &  &  &  & 6/1 & \\
 &  & 8 & -1 &  & -1 &  & -1 &  &  & -2 &  &  &  & 6/1 & \\
 &  & 8 & 3 &  & 3 &  & 3 &  &  & -2 &  &  &  & 6/1 & \\
 &  & 16 & -8 &  & -8 &  & 4 &  &  & -1 &  &  &  & 6/1 & \\
 &  & 16 & -3 &  & -3 &  &  &  &  &  &  &  &  & 6/1 & \\
 &  & 24 &  &  & -9 &  & 6 &  &  &  & -1 &  &  & 6/1 & \\
 &  & 32 & -20 & 3 &  &  &  &  &  &  &  &  &  & 6/1\footnote{[0] Sections 4.4, 4.11.6} & $S$\\
 &  & 32 & -12 &  & -24 & 9 &  &  &  &  &  &  &  & 6/1 & $I$\\
 & 4 &  & -2 &  & -2 & -1 & 1 &  &  &  &  &  &  & 6/1 & \\
 & 4 &  &  &  &  & -1 &  &  &  &  &  &  &  & 6/1 & \\
 & 6 & 4 & -5 & 1 &  &  &  &  &  &  &  &  &  & 6/1 & $S$\\
 & 8 &  & -4 &  &  & 1 &  &  &  &  &  &  &  & 6/1 & $I$\\
 & 16 & -8 &  &  &  & -4 & 4 &  &  & -1 &  &  &  & 6/1 & \\
 & 24 &  & -4 & 1 &  & 3 &  &  &  &  &  &  &  & 6/1 & \\
 & 48 &  &  &  &  & -12 &  & 4 & -1 &  &  &  &  & 6/1 & \\
16 & -12 & -8 &  &  &  & 2 & 3 &  &  & 1 &  &  &  & 6/1 & $I$\\
16 & -8 & -8 & & & & 1 & -2 & & & 1 & & & & 6/1 &\\
16 & -8 &  & 12 &  &  & 1 &  &  &  &  &  &  &  & 6/1 & $I$\\
16 & -1 &  &  &  &  &  &  &  &  &  &  &  &  & 6/1 & $S$\\
16 &  & -8 &  &  &  & -1 &  &  &  & 1 &  &  &  & 6/1 & \\
16 &  &  &  &  &  & -4 & 4 &  &  & -1 &  &  &  & 6/1 & \\
16 & 8 & -8 &  &  &  &  & -2 &  &  & 1 &  &  &  & 6/1 & \\
16 & 8 & -8 & 4 &  & 4 &  & 2 &  &  & 1 &  &  &  & 6/1 & $I$\\
32 & 4 & -8 & 1 &  & 1 &  &  &  &  &  &  &  &  & 6/1 & $I$\\
64 & -20 &  &  &  &  & 1 &  &  &  &  &  &  &  & 6/1 & \\
144 &  & -40 &  &  &  &  &  &  &  & 1 &  &  &  & 6/1 & \\
256 & -64 &  & 16 & -3 &  &  &  &  &  &  &  &  &  & 6/1\footnote{[0] Section 4.2.6, arr. 287} & $S$\\
\hline
1 &  & -2 & 2 &  & 4 &  & -4 &  &  & 1 & 2 & 1 &  & 7/1 & \\
\hline
 &  &  &  &  &  & 2 & -1 &  &  &  &  &  &  & 8/1 & \\
 &  &  &  &  &  & 8 &  & -8 & 2 & -2 & 4 & -1 &  & 8/1 & \\
 &  &  &  &  & 16 & -4 & -16 & 8 & -1 &  &  &  &  & 8/1 & \\
 &  &  & 2 & -1 & -2 & 1 &  &  &  &  &  &  &  & 8/1 & \\
 &  &  & 2 &  & 2 & -1 & -2 &  &  &  &  &  &  & 8/1 & \\
 &  &  & 8 &  & 8 & 1 & -12 &  &  & 4 &  &  &  & 8/1 & \\
 &  & 8 & -6 & 1 & -2 & 1 &  &  &  &  &  &  &  & 8/1 & \\
 &  & 32 & -16 & 2 & -4 & 1 &  &  &  &  &  &  &  & 8/1 & \\
 & 1 &  & -2 & -1 &  & -1 &  & 3 & 1 &  &  & -2 & -1 & 8/1 & $L_2$\\
 & 1 &  & -1 &  & 3 &  & -7 &  &  &  & 4 &  &  & 8/1 & \\
 & 1 &  &  &  &  &  &  &  &  &  &  &  &  & 8/1\footnote{[0] Sections 4.8.1, 4.8.3} & $SL_2$\\
 & 2 &  & -3 &  &  &  &  &  &  &  &  &  &  & 8/1\footnote{[0] Section 4.8.3} & $SI$\\
 & 4 &  & -6 &  & -6 &  & 9 &  &  &  &  &  &  & 8/1 & \\
 & 4 &  & -2 &  & -2 &  & 1 &  &  &  &  &  &  & 8/1 & \\
 & 4 &  & -1 &  & -1 & 1 &  &  &  &  &  &  &  & 8/1 & \\
 & 4 &  &  &  &  & -5 & -5 &  &  &  &  &  &  & 8/1 & \\
 & 4 &  & 3 &  & 3 & -3 &  &  &  &  &  &  &  & 8/1 & \\
 & 8 &  & -4 & 1 &  &  &  &  &  &  &  &  &  & 8/1\footnote{[0] Sections 4.2.3, 4.8.3, arr. 238} & $SC_1$\\
 & 12 &  &  &  &  & 1 & 1 &  &  &  &  &  &  & 8/1 & \\
 & 16 &  & -4 & 1 & 12 & -3 &  &  &  &  &  &  &  & 8/1 & \\
 & 16 &  & -2 & -1 & -6 & 3 &  &  &  &  &  &  &  & 8/1 & \\
 & 24 & -12 &  &  &  & -2 & 1 &  &  &  &  &  &  & 8/1 & \\
4 & -8 & & -2 & 2 & 6 & -3 & & & & & & & & 8/1 & \\
4 &  & -9 & 6 & -1 &  &  &  &  &  &  &  &  &  & 8/1 & $S$\\
4 &  & -1 &  &  &  &  &  &  &  &  &  &  &  & 8/1 & $S$\\
4 & 3 & -1 &  &  &  &  &  &  &  &  &  &  &  & 8/1 & $S$\\
4 & 12 & -1 & 2 & -1 &  &  &  &  &  &  &  &  &  & 8/1 & $S$\\
16 & -12 & 8 &  &  &  &  & 3 &  &  & -3 &  &  &  & 8/1 & $I$\\
16 & -8 & -8 & 6 &  & 6 &  & -4 &  &  & 1 &  &  &  & 8/1 & $I$\\
16 & -8 &  & 4 & -1 &  &  &  &  &  &  &  &  &  & 8/1\footnote{[0] Section 4.2.3, arr. 238} & $SC_1$\\
16 & -6 & -4 & 5 & -1 &  &  &  &  &  &  &  &  &  & 8/1 & $S$\\
16 &  & -3 & 1 &  &  &  &  &  &  &  &  &  &  & 8/1 & $SI$\\
16 &  & 12 & -7 & 1 &  &  &  &  &  &  &  &  &  & 8/1 & $S$\\
16 & 4 & -8 &  &  &  &  & -1 &  &  & 1 &  &  &  & 8/1 & \\
16 & 4 & -8 & 2 &  & 2 &  & 1 &  &  & 1 &  &  &  & 8/1 & \\
16 & 8 &  & -2 &  & -2 & 1 &  &  &  &  &  &  &  & 8/1 & \\
32 & -4 & -8 & 3 &  & 3 &  &  &  &  &  &  &  &  & 8/1 & $I$\\
32 & -4 & 1 &  &  &  &  &  &  &  &  &  &  &  & 8/1 & $S$\\
32 &  &  & -3 & 1 &  &  &  &  &  &  &  &  &  & 8/1 & $S$\\
64 & -12 & -16 &  &  &  &  & 3 &  &  &  &  &  &  & 8/1 & \\
64 &  & -16 & 8 & -1 &  &  &  &  &  &  &  &  &  & 8/1 & $S$\\
256 &  &  &  & -1 &  &  &  &  &  &  &  &  &  & 8/1\footnote{[0] Section 4.7.2} & $S$\\
\hline
 &  &  &  &  & 12 &  & -12 & 3 &  &  & 3 & -1 &  & 9/1 & \\
 &  &  &  &  & 32 &  &  & 12 &  &  & -8 & -2 & 1 & 9/1 & \\
 &  & 4 & -2 &  & -2 & -3 & 4 &  &  & -1 &  &  &  & 9/1 & \\
 &  & 64 & -16 &  & -16 & 3 &  &  &  &  &  &  &  & 9/1 & \\
 & 4 &  & 2 &  & 2 & -1 & 1 &  &  &  &  &  &  & 9/1 & \\
 & 12 & -4 &  &  &  & 9 & -6 &  &  & 1 &  &  &  & 9/1 & \\
16 & -16 & 4 & 2 &  & 2 & 3 & -4 &  &  & 1 &  &  &  & 9/1 & $I$\\
16 & -12 & -4 & 1 & & 1 & 2 & & & & & & & & 9/1 &\\
16 & -4 & -8 & 2 &  & 2 &  & -1 &  &  & 1 &  &  &  & 9/1 & \\
16 & -4 & -4 & -1 &  & -1 &  &  &  &  &  &  &  &  & 9/1 & \\
16 & & -16 & 4 & & 4 & -1 & & & & & & & & 9/1 &\\
16 & 4 & -8 & 2 & & 2 & 1 & -3 & & & 1 & & & & 9/1 &\\
16 & 4 & -4 &  &  &  & -2 & 1 &  &  &  &  &  &  & 9/1 & \\
32 & -8 & -2 & 1 &  & 1 &  &  &  &  &  &  &  &  & 9/1 & $I$\\
32 &  &  &  & -4 &  & -2 &  &  & 1 &  &  &  &  & 9/1 & \\
256 & 32 &  &  & -3 &  &  &  &  &  &  &  &  &  & 9/1 & $S$\\
\hline
 &  &  & 4 &  & -4 & 1 &  & -2 &  &  &  & 1 &  & 10/1 & \\
 & 1 &  & -1 &  &  &  &  &  &  &  &  &  &  & 10/1\footnote{[0] Sections 4.8.3, 4.9} & $SI$\\
 & 1 &  & 3 &  &  &  & -4 &  &  &  &  &  &  & 10/1 & $I$\\
16 & 8 & -4 &  &  &  & 1 &  &  &  &  &  &  &  & 10/1 & \\
64 & -18 &  & 5 & -1 &  &  &  &  &  &  &  &  &  & 10/1 & $S$\\
256 &  &  & -16 & 5 &  &  &  &  &  &  &  &  &  & 10/1 & $S$\\
\hline
 &  &  &  &  &  &  & 32 & -20 & 3 &  &  &  &  & 12/1 & \\
 &  &  &  &  & 8 &  & -8 & 2 &  &  & -6 & -3 &  & 12/1 & \\
 &  &  & 1 &  & 1 &  & 1 &  &  &  &  &  &  & 12/1 & \\
 &  &  & 2 &  & 2 &  & -6 &  &  &  & 4 & -1 &  & 12/1 & \\
 &  & 4 & -1 &  & -1 &  &  &  &  &  &  &  &  & 12/1 & \\
 &  & 4 & 4 & 1 &  & -3 &  &  &  &  &  &  &  & 12/1 & \\
 &  & 64 &  &  &  & -27 & 8 &  &  &  &  &  &  & 12/1 & \\
 & 4 &  & 2 &  & 2 & -4 & 1 &  &  &  &  &  &  & 12/1 & \\
 & 8 & -4 &  &  &  & 4 & -4 &  &  & 1 &  &  &  & 12/1 & \\
 & 16 &  & -8 &  & -8 & -1 & 4 &  &  &  &  &  &  & 12/1 & \\
 & 16 &  &  &  &  & -1 &  &  &  &  &  &  &  & 12/1 & \\
 & 32 &  & -16 & 3 &  & 1 &  &  &  &  &  &  &  & 12/1 & \\
1 & -1 &  & 1 &  &  &  &  &  &  &  &  &  &  & 12/1\footnote{[0] Section 4.2.3, arr. 239} & $SIL_1$\\
4 & -1 &  &  &  &  &  &  &  &  &  &  &  &  & 12/1 & $S$\\
12 & -6 & -4 & 5 & -1 &  &  &  &  &  &  &  &  &  & 12/1 & $S$\\
12 & 12 & -4 &  &  &  & 3 & -2 &  &  &  &  &  &  & 12/1 & \\
16 & -8 & -4 &  &  &  & 1 & 1 &  &  &  &  &  &  & 12/1 & \\
16 & -4 & -8 & -2 &  & -2 &  & -1 &  &  & 1 &  &  &  & 12/1 & \\
16 & -4 & -8 &  &  &  &  & 1 &  &  & 1 &  &  &  & 12/1 & \\
16 & -4 & -4 & 1 &  & 1 &  &  &  &  &  &  &  &  & 12/1 & \\
16 & 12 & -8 &  &  &  & 2 & -3 &  &  & 1 &  &  &  & 12/1 & $I$\\
16 & 28 & -4 &  &  &  & -8 & 1 &  &  &  &  &  &  & 12/1 & $I$\\
48 & -12 &  & 4 & -1 &  &  &  &  &  &  &  &  &  & 12/1 & $S$\\
64 & 16 & -4 &  &  &  & 1 &  &  &  &  &  &  &  & 12/1 & \\
81 & -27 &  & 9 & -2 &  &  &  &  &  &  &  &  &  & 12/1 & $SL_1$\\
144 &  &  & -4 & 1 &  &  &  &  &  &  &  &  &  & 12/1 & $S$\\
\hline
 &  &  &  &  &  & 2 & -9 &  &  & 4 &  &  &  & 14/1 & \\
 & 8 & -4 & 2 &  & 2 & -9 & 2 &  &  &  &  &  &  & 14/1 & \\
 & 32 &  & 12 &  &  & -7 & -14 &  &  &  &  &  &  & 14/1 & $I$\\
16 &  & -8 &  &  &  &  &  &  &  & 1 & 8 &  &  & 14/1 & \\
64 & 12 &  &  &  &  & 1 &  &  &  &  &  &  &  & 14/1 & \\
\hline
 &  & 16 &  &  &  & -1 & -4 &  &  & -4 &  &  &  & 14/2 & \\
 & 2 &  & -1 &  &  &  &  &  &  &  &  &  &  & 14/2\footnote{[0] Section 4.8.3} & $SI$\\
 & 4 & -8 & 3 &  & 3 &  & 2 &  &  & 2 &  &  &  & 14/2 & \\
128 &  &  & -4 & 1 &  &  &  &  &  &  &  &  &  & 14/2 & $S$\\
\hline
1 & -4 & -2 & 4 &  & 4 &  & -4 &  &  & 1 &  &  &  & 15/1 & \\
45 & -10 &  & -4 &  &  & 1 &  &  &  &  &  &  &  & 15/1 & $I$\\
\hline
1 & -2 & -4 &  &  &  & 1 &  &  &  &  &  &  &  & 15/2 & \\
9 & -18 & 4 &  &  &  & 9 & -4 &  &  &  &  &  &  & 15/2 & $I$\\
9 &  &  & -4 &  &  &  &  &  &  &  &  &  &  & 15/2 & $SI$\\
81 &  & -8 &  &  &  &  & 8 &  &  &  & 4 &  &  & 15/2 & \\
225 &  & 16 & 60 &  &  &  &  &  &  &  &  &  &  & 15/2 & $SI$\\
\hline
 &  & 4 & -4 & 1 &  & -4 & 4 &  &  & -1 &  &  &  & 20/1 & \\
 &  & 20 & -9 &  & -9 & 4 &  &  &  &  &  &  &  & 20/1 & \\
 & 4 & -16 & -6 & & -6 & & 5 & & & 4 & & & & 20/1 &\\
 & 4 &  & -8 &  &  &  & -1 &  &  &  &  &  &  & 20/1 & $I$\\
 & 16 & -12 & 3 &  & 3 &  &  &  &  &  &  &  &  & 20/1 & \\
1 & -7 &  & 5 &  &  & 1 &  &  &  &  &  &  &  & 20/1 & $I$\\
3 & 1 & -12 & 3 &  &  &  &  &  &  &  &  &  &  & 20/1 & $SI$\\
4 & -2 & 4 & -1 &  &  &  &  &  &  &  &  &  &  & 20/1 & $SI$\\
5 & -3 & -20 & 13 & -2 &  &  &  &  &  &  &  &  &  & 20/1 & $S$\\
16 &  & -4 &  &  &  &  & 3 &  &  &  &  &  &  & 20/1 & \\
20 & -9 &  & 4 & -1 &  &  &  &  &  &  &  &  &  & 20/1 & $S$\\
\hline
 &  & 4 &  &  &  & -9 & 6 &  &  & -1 &  &  &  & 21/1 & \\
 &  & 4 & 12 &  &  &  &  &  &  & -1 &  &  &  & 21/1 & $I$\\
 &  & 64 &  &  &  & -9 &  &  &  &  &  &  &  & 21/1 & \\
 & 72 & -36 &  &  &  &  & -2 &  &  & 1 &  &  &  & 21/1 & $I$\\
16 & 4 & -8 & & & & 7 & -5 & & & 1 & & & & 21/1 & \\
48 & -12 & -12 &  &  &  &  & 1 &  &  &  &  &  &  & 21/1 & \\
\hline
9 & -6 &  & 16 & -4 &  & -3 &  &  &  &  &  &  &  & 21/2 & \\
25 & 20 & -4 & -4 &  & -4 &  &  &  &  &  &  &  &  & 21/2 & $I$\\
27 &  & 18 &  &  &  &  &  &  &  & 3 & -4 &  &  & 21/2 & \\
\hline
16 &  & -8 &  &  &  &  & 16 &  &  & 1 & -16 &  &  & 22/1 & \\
\hline
 &  &  &  &  & 128 & -64 &  &  & 1 &  &  &  &  & 22/2 & \\
 &  &  & 4 & -1 & 4 & 1 & -4 & 1 &  &  &  &  &  & 22/2 & \\
 & 64 & -64 & 16 & -1 &  &  &  &  &  &  &  &  &  & 22/2 & $S$\\
\hline
 &  &  &  &  & 64 & -16 &  & -4 & 1 &  &  &  &  & 23/1 & \\
 & 16 &  &  & -1 &  & -5 &  & 1 &  &  &  &  &  & 23/1 & \\
\hline
 &  &  &  &  &  &  & 16 & -6 &  &  & -24 & 9 &  & 24/1 & \\
 &  &  &  &  &  &  & 256 & -168 & 27 &  &  &  &  & 24/1 & \\
 &  &  &  &  &  & 9 & -4 &  &  &  &  &  &  & 24/1 & \\
 &  &  &  &  & 8 &  & -8 & 2 &  &  & 2 & -1 &  & 24/1 & \\
 &  &  &  &  & 12 &  & -20 &  &  &  & 3 &  &  & 24/1 & \\
 &  &  &  &  & 24 &  & 8 & 6 &  &  & -2 & -1 &  & 24/1 & \\
 &  &  & 3 &  & 3 &  & -1 &  &  &  &  &  &  & 24/1 & \\
 &  &  & 6 &  & 6 &  & -10 &  &  &  & 4 & -1 &  & 24/1 & \\
 &  &  & 8 &  &  &  & 16 & -8 &  &  & -2 & 1 &  & 24/1 & \\
 &  &  & 16 &  &  &  & -16 & 2 &  &  & 4 & -1 &  & 24/1 & \\
 &  & 1 & 3 &  &  &  & -1 & -3 &  &  & 1 & 3 &  & 24/1 & $L_3$\\
 &  & 3 & -1 &  &  &  &  &  &  &  &  &  &  & 24/1\footnote{[0] Sections 4.4, 4.11.6} & $SIL_3$\\
 &  & 4 & -3 &  & -3 &  & 3 &  &  & -1 &  &  &  & 24/1 & \\
 &  & 4 & 1 &  & -1 &  &  &  &  &  &  &  &  & 24/1 & $I$\\
 &  & 64 & -24 &  & -24 & 9 &  &  &  &  &  &  &  & 24/1 & \\
 & 4 &  &  & -1 & 4 & -4 &  & -1 & 1 &  &  &  &  & 24/1 & \\
 & 8 & & -8 & 4 & & -2 & & 2 & -1 & & & & & 24/1 & \\
 & 12 &  &  &  &  & -3 & 1 &  &  &  &  &  &  & 24/1 & \\
 & 16 &  &  & -4 &  & -4 &  &  & 1 &  &  &  &  & 24/1 & \\
 & 24 & -12 &  &  &  & -6 & 5 &  &  & -1 &  &  &  & 24/1 & $I$\\
 & 36 &  & -6 &  & -6 &  & 1 &  &  &  &  &  &  & 24/1 & \\
 & 96 & -64 &  &  &  & 3 &  &  &  &  &  &  &  & 24/1 & \\
4 & -8 & & 5 & & & 4 & & -5 & & & & 1 & & 24/1 & \\
12 & -3 & 1 &  &  &  &  &  &  &  &  &  &  &  & 24/1 & $S$\\
16 & -8 & -8 & 6 &  & 6 & 1 & -4 &  &  & 1 &  &  &  & 24/1 & $I$\\
16 & -8 & -4 &  &  &  & 1 &  &  &  &  &  &  &  & 24/1 & \\
16 & -8 & 8 & -6 &  & -6 & 1 &  &  &  & 1 &  &  &  & 24/1 & $I$\\
16 & -8 & 12 & -3 &  & -3 & 1 &  &  &  &  &  &  &  & 24/1 & $I$\\
16 & 4 & -8 &  &  &  & -2 & -1 &  &  & 1 &  &  &  & 24/1 & \\
32 & -4 & 1 & 6 &  &  & -1 &  &  &  &  &  &  &  & 24/1 & $I$\\
36 &  & -9 & 6 & -1 &  &  &  &  &  &  &  &  &  & 24/1 & $S$\\
36 &  & -1 &  &  &  &  &  &  &  &  &  &  &  & 24/1 & $S$\\
36 &  &  & -3 & 1 &  &  &  &  &  &  &  &  &  & 24/1 & $S$\\
36 &  & 16 & -3 &  &  &  &  &  &  &  &  &  &  & 24/1 & $SI$\\
64 & 16 & -16 &  &  &  & 1 &  &  &  &  &  &  &  & 24/1 & \\
192 &  & 16 & -8 & 1 &  &  &  &  &  &  &  &  &  & 24/1 & $S$\\
\hline
 &  & 32 & -12 & 1 & -8 & 2 &  &  &  &  &  &  &  & 26/2 & \\
 & 2 & -16 & 8 & -1 & 3 & -1 &  &  &  &  &  &  &  & 26/2 & \\
 & 8 &  & -4 & 1 &  & -1 &  &  &  &  &  &  &  & 26/2 & \\
\hline
16 & & -8 & & 1 & & 1 & & -1 & & 1 & & & & 28/1 & \\
\hline
 &  & 4 &  &  &  & -1 &  &  &  &  &  &  &  & 28/2 & \\
& 16 & -4 & -4 & & & 9 & -6 & & & 1 & & & & 28/2 & $I$\\
1 &  &  &  &  & 1 &  & -2 &  &  &  & 1 &  &  & 28/2 & \\
16 &  &  & -4 & 1 &  &  &  &  &  &  &  &  &  & 28/2 & $S$\\
16 & 12 & -1 & 2 & -1 &  &  &  &  &  &  &  &  &  & 28/2 & $S$\\
49 & -9 &  & 1 &  &  &  &  &  &  &  &  &  &  & 28/2 & $SI$\\
\hline
 &  & 16 & 5 &  & 5 &  &  &  &  &  &  &  &  & 30/1 & \\
16 & 56 &  & -20 &  &  & 9 &  &  &  &  &  &  &  & 30/1 & $I$\\
64 & -52 &  &  &  &  & 9 &  &  &  &  &  &  &  & 30/1 & \\
80 & 1 &  & -1 &  &  &  &  &  &  &  &  &  &  & 30/1 & $SI$\\
\hline
 &  & 32 & -18 &  &  &  & 1 &  &  &  &  &  &  & 30/2 & $I$\\
 &  & 256 &  &  & -64 &  &  &  &  &  &  & -27 &  & 30/2 & \\
 & 2 &  & -2 &  & -1 & 1 & 2 & -2 &  &  & -1 & 1 &  & 30/2 & \\
 & 32 &  & -4 &  &  & -5 & 2 &  &  &  &  &  &  & 30/2 & $I$\\
 & 36 & -16 &  &  &  & -9 & 4 &  &  &  &  &  &  & 30/2 & \\
 & 81 & -36 &  &  &  &  & -9 &  &  & 4 &  &  &  & 30/2 & $I$\\
64 & -10 &  & 1 &  &  &  &  &  &  &  &  &  &  & 30/2 & $SI$\\
\hline
 &  &  & 2 &  &  &  & -1 & 1 &  &  & 1 & -1 &  & 32/1 & \\
 &  &  & 2 & 1 &  &  & 1 & -1 &  &  &  &  &  & 32/1 & \\
 &  &  & 4 & 5 &  &  & -2 & -3 & -1 &  &  &  &  & 32/1 & \\
 &  &  & 16 &  &  &  & 8 & -16 &  &  &  & 6 & -3 & 32/1 & \\
 &  &  & 27 &  &  &  &  &  &  & -4 &  &  &  & 32/1 & $I$\\
\hline
 &  &  & 4 & -1 & 4 & -1 & -4 & 1 &  &  &  &  &  & 32/2\footnote{[0] Sections 4.3.1, B.10} & \\
 &  &  & 8 & -1 & 8 & -1 & -8 & 1 & & & & & & 32/2 & \\
 &  & 16 & -8 & 1 & -4 & 1 &  &  &  &  &  &  &  & 32/2 & \\
 &  & 16 & -4 &  & -4 & 1 &  &  &  &  &  &  &  & 32/2 & \\
 & 2 & -1 &  &  &  &  &  &  &  &  &  &  &  & 32/2 & $SC_2$\\
 & 4 & -16 &  &  &  & -2 & 7 &  &  & 4 &  &  &  & 32/2 & \\
 & 4 & -4 & 1 &  & -3 & 1 &  &  &  &  &  &  &  & 32/2 & $I$\\
 & 4 &  &  &  &  &  & -1 &  &  &  &  &  &  & 32/2 & \\
 & 4 &  &  &  &  & 1 & -1 &  &  &  &  &  &  & 32/2 & \\
 & 4 &  & 2 &  & -6 &  & -3 &  &  &  &  &  &  & 32/2 & $I$\\
 & 4 &  & 2 &  & 2 &  & 1 &  &  &  &  &  &  & 32/2 & \\
 & 8 & -16 & 6 &  & 6 & -3 &  &  &  &  &  &  &  & 32/2 & \\
 & 8 &  & -6 &  & -6 & -1 & 4 &  &  &  &  &  &  & 32/2 & \\
 & 8 &  & -3 &  & -3 &  & 1 &  &  &  &  &  &  & 32/2 & \\
4 & -6 & -1 & 4 & -1 &  &  &  &  &  &  &  &  &  & 32/2 & $SC_3$\\
4 & -2 & 1 &  &  &  &  &  &  &  &  &  &  &  & 32/2 & $SC_2$\\
8 & -6 & -1 & 4 & -1 &  &  &  &  &  &  &  &  &  & 32/2 & $SC_3$\\
8 &  & 6 & -5 & 1 &  &  &  &  &  &  &  &  &  & 32/2 & $S$\\
16 & -8 & -8 & 6 &  & 6 & -3 &  &  &  &  &  &  &  & 32/2 & $I$\\
16 &  & -4 & 4 & -1 &  &  &  &  &  &  &  &  &  & 32/2\footnote{[0] Section 4.7.1} & $S$\\
16 &  & 4 & -4 & 1 &  &  &  &  &  &  &  &  &  & 32/2\footnote{[0] Section 4.7.1} & $S$\\
\hline
1 & -2 & -2 &  &  &  & 1 & -2 &  &  & 1 &  &  &  & 33/1 & \\
\hline
27 &  & -4 &  &  &  &  & 4 &  &  &  &  &  &  & 33/2 & \\
\hline
 &  & 8 & -9 &  & -9 & 4 & 7 & & & -2 & & & & 34/2 & \\
 &  & 64 & -48 & 17 &  & -8 &  &  &  &  &  &  &  & 34/2 & \\
 & 2 & -4 & 1 &  &  &  &  &  &  &  &  &  &  & 34/2 & $SI$\\
 & 4 &  & -8 &  &  & 1 & -1 &  &  &  &  &  &  & 34/2 & $I$\\
 & 8 &  & -12 &  &  & -4 &  & 8 &  &  &  & -3 &  & 34/2 & \\
16 & -9 &  & 4 & -1 &  &  &  &  &  &  &  &  &  & 34/2 & $S$\\
\hline
 &  & 64 & -48 & 8 &  & 1 &  &  &  &  &  &  &  & 38/1 & \\
\hline
 & 4 & 4 &  &  &  & -1 & -2 &  &  & -1 &  &  &  & 39/1 & \\
16 & 12 & -8 & 8 &  &  & -1 & -3 &  &  & 1 &  &  &  & 39/1 & $I$\\
\hline
16 & -4 & -4 &  &  &  & -1 & 3 &  &  & -1 &  &  &  & 40/1 & \\
\hline
 &  &  &  &  &  & 4 & -8 & -12 &  &  & 12 & 9 &  & 40/3 & \\
 &  &  & 1 &  & 1 & -1 &  &  &  &  &  &  &  & 40/3 & \\
 &  &  & 16 &  &  & -1 & 4 &  &  & -4 &  &  &  & 40/3 & $I$\\
 & 4 & -4 &  &  &  &  & -1 &  &  & 1 &  &  &  & 40/3 & \\
 & 4 &  & -4 & 1 &  &  &  &  &  &  &  &  &  & 40/3\footnote{[0] Sections 4.8.1, 4.8.3} & $S$\\
 & 8 &  &  & -1 & -4 & 2 &  &  &  &  &  &  &  & 40/3 & \\
 & 16 &  &  &  &  & -8 &  & 4 & -1 &  &  &  &  & 40/3 & \\
 & 24 & -12 &  &  &  & 10 & -11 &  &  & 3 &  &  &  & 40/3 & $I$\\
4 &  &  & 1 &  &  &  &  &  &  &  &  &  &  & 40/3 & $SI$\\
4 & 4 & -1 &  &  &  &  &  &  &  &  &  &  &  & 40/3 & $S$\\
64 & -48 & -16 &  &  &  & 9 &  &  &  &  &  &  &  & 40/3 & $I$\\
64 & -32 & -16 & 24 & -5 &  &  &  &  &  &  &  &  &  & 40/3 & $S$\\
\hline
 & 7 &  & 1 &  &  & -16 &  &  &  &  &  &  &  & 42/1 & $I$\\
\hline
 & 2 &  & -1 &  &  & -2 &  & 3 &  &  &  & -1 &  & 42/2 & \\
 & 3 &  & -5 &  &  & -3 &  & 8 &  &  &  & -5 &  & 42/2 & \\
 & 9 & -4 &  &  &  &  &  &  &  &  &  &  &  & 42/2 & $S$\\
 & 9 &  & -1 &  &  &  &  &  &  &  &  &  &  & 42/2\footnote{[0] Section 4.8.3} & $SI$\\
 & 36 &  & 6 &  & 6 & -9 & 1 &  &  &  &  &  &  & 42/2 & $I$\\
 & 54 &  & -27 & 7 &  &  &  &  &  &  &  &  &  & 42/2\footnote{[0] Section 4.8.3} & $S$\\
\hline
 & 4 & 16 & -6 &  & -6 &  & 1 &  &  &  &  &  &  & 44/1 & \\
\hline
3 &  & -12 &  &  &  &  & 4 &  &  &  &  &  &  & 45/1 & \\
\hline
9 & -6 &  & 10 &  &  & -3 &  & 2 &  &  &  & 1 &  & 51/2 & \\
17 & -8 & 16 & -4 &  &  &  &  &  &  &  &  &  &  & 51/2 & $SI$\\
81 & -36 &  & 16 & -4 &  &  &  &  &  &  &  &  &  & 51/2 & $S$\\
\hline
 &  & 64 & -32 & 4 &  & -1 &  &  &  &  &  &  &  & 53/1 & \\
\hline
 &  &  &  &  &  &  &  64 & -6 & & & -32 & 9 & & 54/1 & \\
 &  &  &  &  & 24 &  & -16 & 6 &  &  & 2 & -1 &  & 54/1 & \\
 &  &  & 4 &  &  & 3 & -2 &  &  &  &  &  &  & 54/1 & $I$\\
 &  & 8 & -1 &  & -2 &  &  &  &  &  &  &  &  & 54/1 & $I$\\
 &  & 32 & -8 &  & -12 & 3 &  &  &  &  &  &  &  & 54/1 & $I$\\
 &  & 64 & -16 &  & -32 & 9 &  &  &  &  &  &  &  & 54/1 & $I$\\
\hline
 &  &  &  &  &  & 9 & -2 & -3 &  &  & 1 &  &  & 54/2 & \\
 &  &  &  &  & 12 &  & -20 & 3 &  &  & 7 & -2 &  & 54/2 & \\
 &  &  &  &  & 12 &  & -4 & 3 &  &  & -1 &  &  & 54/2 & \\
 &  &  &  &  & 128 &  &  & -32 & 9 &  &  &  &  & 54/2 & \\
 &  & 64 & -16 &  &  & -3 &  &  &  &  &  &  &  & 54/2 & $I$\\
16 & -8 & -8 & 8 &  & 4 & 1 & -4 &  &  & 1 &  &  &  & 54/2 & $I$\\
\hline
 &  &  &  &  &  & 4 & -16 & 8 & -1 &  &  &  &  & 56/2 & \\
 & 4 & -1 &  &  &  &  &  &  &  &  &  &  &  & 56/2 & $S$\\
 & 64 &  &  &  &  &  & -16 & 8 & -1 &  &  &  &  & 56/2 & \\
4 & -8 & -12 & 8 &  & 8 &  & -4 &  &  & 1 &  &  &  & 56/2 & $I$\\
\hline
27 &  &  &  &  &  &  &  &  &  &  & 4 &  &  & 57/1 & \\
\hline
16 & 3 & -4 & 1 &  &  &  &  &  &  &  &  &  &  & 58/2 & $SI$\\
32 &  &  & -2 & 1 & 2 & -1 &  &  &  &  &  &  &  & 58/2 & \\
\hline
 &  &  &  &  & 96 & 9 & -56 &  &  & 16 &  &  &  & 60/1 & $I$\\
&  &  & 32 &  & 32 & -27 &  &  &  &  &  &  &  & 60/1 & \\
& 8 & & -4 & & & 10 & -2 & -3 & & & 1 & & & 60/1 & \\
 & 80 & -64 &  &  &  & -5 & 4 &  &  &  &  &  &  & 60/1 & $I$\\
16 & -8 & 16 & -4 &  & -4 & 1 &  &  &  &  &  &  &  & 60/1 & $I$\\
16 & -4 & -4 &  &  &  &  & 1 &  &  &  &  &  &  & 60/1 & \\
16 & -4 & 1 &  &  &  &  &  &  &  &  &  &  &  & 60/1 & $S$\\
\hline
& 16 &  & 64 &  &  & -25 &  &  &  &  &  &  &  & 60/2 & $I$\\
36 & -9 & -4 &  &  &  &  & 1 &  &  &  &  &  &  & 60/2 & \\
324 & 27 &  &  & -1 &  &  &  &  &  &  &  &  &  & 60/2 & $S$\\
432 &  & -28 &  &  &  &  & 1 &  &  &  &  &  &  & 60/2 & $I$\\
\hline
16 &  & -8 &  &  &  &  & 8 &  &  & 1 & -4 &  &  & 62/1 & \\
\hline
&  & 4 & 4 &  &  & -1 & -2 &  &  & -1 &  &  &  & 68/1 & $I$\\
\hline
&  &  &  &  & 16 &  & -16 & 2 &  &  & 4 & -1 &  & 72/1 & \\
 &  & 4 & -12 & 3 &  & 6 &  &  &  &  &  &  &  & 72/1 & \\
 &  & 16 &  &  &  & -3 &  &  &  &  &  &  &  & 72/1 & \\
 & 3 & -1 &  &  &  &  &  &  &  &  &  &  &  & 72/1 & $S$\\
 & 12 &  & -2 &  & -18 &  & 3 &  &  &  &  &  &  & 72/1 & $I$\\
 & 12 &  &  &  &  &  & 1 &  &  &  &  &  &  & 72/1 & \\
12 &  & -3 & 1 &  &  &  &  &  &  &  &  &  &  & 72/1 & $SI$\\
\hline
1 &  & 14 &  &  &  &  &  &  &  & 1 & 4 &  &  & 73/1 & \\
\hline
 &  &  & 4 & -3 & 4 & -1 & -4 & 7 & -2 &  &  &  &  & 74/2 & \\
&  & 64 & -16 &  &  & 1 &  &  &  &  &  &  &  & 74/2 & $I$\\
16 &  & -8 & 4 &  &  &  & -2 &  &  & 1 &  &  &  & 74/2 & $I$\\
\hline
& 108 & -64 &  &  &  &  & 1 &  &  &  &  &  &  & 78/2 & \\
\hline
32 &  & -8 & -22 &  &  &  & -3 &  &  &  &  &  &  & 88/1 & $I$\\
\hline
& 4 &  & 1 &  & 1 &  &  &  &  &  &  &  &  & 88/2 & \\
 & 16 &  & -4 & 1 & -4 & 1 &  &  &  &  &  &  &  & 88/2 & \\
 & 16 &  &  & -1 &  &  &  &  &  &  &  &  &  & 88/2\footnote{[0] Section 4.8.3} & $S$\\
 & 16 &  & 4 & -1 & -12 & -5 &  &  &  &  &  &  &  & 88/2 & \\
 & 256 &  &  &  &  & -32 &  &  & 1 &  &  &  &  & 88/2 & \\
\hline
&  & 32 & -8 &  & -10 &  & 1 &  &  &  &  &  &  & 90/1 & $I$\\
 & 3 & 4 & -1 &  & -1 &  &  &  &  &  &  &  &  & 90/1 & \\
 & 12 & 8 & 5 &  & 5 &  & 2 &  &  & -2 &  &  &  & 90/1 & \\
\hline
1 &  &  & 2 &  & -4 &  &  &  &  &  &  & 1 &  & 95/4 & \\
1 & 2 &  & -4 &  &  & 1 &  &  &  &  &  &  &  & 95/4 & $I$\\
\hline
 &  &  &  &  & 27 &  & -54 &  &  & 4 & 27 &  &  & 96/1 & \\
&  &  & 2 &  &  &  & -3 &  &  &  &  &  &  & 96/1 & $I$\\
 &  &  & 2 &  &  &  & 1 & -2 &  &  & 2 & 4 &  & 96/1 & \\
 &  &  & 2 &  & 6 &  & -6 & 1 &  &  & -3 & -1 &  & 96/1 & \\
 &  &  & 4 & -1 &  &  & 6 & -5 & 1 &  &  &  &  & 96/1 & \\
 &  &  & 4 & 3 &  &  & -2 & -1 & -3 &  & 4 &  &  & 96/1 & \\
 &  &  & 6 &  &  &  & -1 &  &  &  &  &  &  & 96/1 & $I$\\
 &  &  & 6 &  & 6 & -9 & 2 &  &  &  &  &  &  & 96/1 & \\
 &  & 8 & -3 &  &  &  & -8 & 3 &  &  & 8 & -3 &  & 96/1 & \\
 &  & 8 & -3 &  &  &  &  &  &  &  &  &  &  & 96/1\footnote{[0] Sections 4.4, 4.11.6} & $SI$\\
 &  & 16 &  &  &  & -9 & 12 &  &  & -4 &  &  &  & 96/1 & \\
 &  & 24 & -17 & 3 &  &  &  &  &  &  &  &  &  & 96/1\footnote{[0] Sections 4.4, 4.11.6} & $S$\\
 & 12 & -2 &  &  &  & -3 &  &  &  &  &  &  &  & 96/1 & \\
 & 36 &  &  &  &  &  & -1 &  &  &  &  &  &  & 96/1 & \\
2 & -4 & 1 &  &  &  & 2 &  &  &  &  &  &  &  & 96/1 & \\
6 & 12 & -1 &  &  &  &  &  &  &  &  &  &  &  & 96/1 & $S$\\
8 & 20 & -2 &  &  &  & -1 &  &  &  &  &  &  &  & 96/1 & \\
16 & -8 & -8 & -3 & & -3 & 1 & -1 & & & 1 & & & & 96/1 & \\
16 & -8 &  & 6 &  & 6 & 1 &  &  &  &  &  &  &  & 96/1 & \\
144 &  & -36 & 20 & -3 &  &  &  &  &  &  &  &  &  & 96/1 & $S$\\
\hline
&  &  &  &  & 4 &  & -8 &  &  &  & 3 &  &  & 96/2 & \\
 &  &  &  &  & 16 &  &  & 8 &  &  & -4 &  & 1 & 96/2 & \\
 &  & 4 & -1 &  &  &  &  &  &  &  &  &  &  & 96/2\footnote{[0] Sections 4.4, 4.11.6} & $SI$\\
 &  & 16 &  &  &  & -27 & 20 &  &  & -4 &  &  &  & 96/2 & \\
 &  & 32 & -9 &  & -9 &  & 1 &  &  &  &  &  &  & 96/2 & \\
24 & 24 & -2 &  &  &  & 6 & -1 &  &  &  &  &  &  & 96/2 & $I$\\
216 &  & 18 & -9 & 1 &  &  &  &  &  &  &  &  &  & 96/2 & $S$\\
\hline
 &  &  &  &  & 8 &  & -24 & 2 &  &  & 18 & -3 &  & 96/4 & \\
&  &  &  &  & 8 &  & 8 & 2 &  &  & 2 & 1 &  & 96/4 & \\
 &  &  &  &  & 12 &  & -8 &  &  &  & 1 &  &  & 96/4 & \\
 &  &  & 1 &  & 1 &  & -1 &  &  &  &  &  &  & 96/4\footnote{[0] Sections 4.3.1, B.10} & \\
 &  &  & 2 &  & 6 &  & -6 & 1 &  &  & 1 & -1 &  & 96/4 & \\
 &  &  & 8 & -3 & -8 & 3 &  &  &  &  &  &  &  & 96/4 & \\
 &  &  & 8 & -3 & 8 & -3 & -8 & 3 &  &  &  &  &  & 96/4 & \\
 &  &  & 16 &  &  &  &  & -4 &  & -8 & 7 & -1 &  & 96/4 & \\
 &  &  & 27 &  &  &  &  &  &  & -2 &  &  &  & 96/4 & $I$\\
 &  &  & 27 &  & 27 &  & -27 &  &  & 8 &  &  &  & 96/4 & \\
 &  & 4 & 3 &  &  &  & -4 & -3 &  &  & 4 & 3 &  & 96/4 & \\
 &  & 12 & -7 & 1 &  &  &  &  &  &  &  &  &  & 96/4\footnote{[0] Sections 4.4, 4.11.6} & $S$\\
 & 4 &  &  &  &  & -1 & -1 &  &  &  &  &  &  & 96/4 & \\
 & 4 &  &  &  &  & 2 & -1 &  &  &  &  &  &  & 96/4 & \\
 & 4 &  & 1 &  & 1 & -1 &  &  &  &  &  &  &  & 96/4 & \\
 & 6 & 1 & -4 & 1 &  &  &  &  &  &  &  &  &  & 96/4 & $SC_4$\\
 & 8 &  & -2 &  & -2 & 1 &  &  &  &  &  &  &  & 96/4 & \\
 & 12 &  & 6 &  & -2 &  & -1 &  &  &  &  &  &  & 96/4 & $I$\\
4 & 2 & -1 &  &  &  &  &  &  &  &  &  &  &  & 96/4 & $SC_5$\\
8 & -2 & 1 &  &  &  &  &  &  &  &  &  &  &  & 96/4 & $SC_5$\\
8 & 8 & -2 &  &  &  & 2 & -1 &  &  &  &  &  &  & 96/4 & \\
12 & -6 & -1 & 4 & -1 &  &  &  &  &  &  &  &  &  & 96/4 & $SC_4$\\
16 & -8 & 8 & -2 & & -2 & 1 & & & & & & & & 96/4 & \\
16 & -4 & -8 & 4 & & 4 & & -3 & & & 1 & & & & 96/4 & \\
16 & -4 &  &  &  &  &  & 1 &  &  & -1 &  &  &  & 96/4 & \\
24 &  & 2 & -1 &  &  &  &  &  &  &  &  &  &  & 96/4 & $SI$\\
64 & -16 & -16 & 16 & -3 &  &  &  &  &  &  &  &  &  & 96/4 & $S$\\
192 & -48 & 16 &  & -1 &  &  &  &  &  &  &  &  &  & 96/4 & $S$\\
\hline
&  & 64 & -48 &  &  & 9 &  &  &  &  &  &  &  & 102/2 & $I$\\
16 & -8 & -8 & 12 &  &  & 1 & -4 &  &  & 1 &  &  &  & 102/2 & $I$\\
\hline
& 4 & -16 & 4 &  & 4 & -1 &  &  &  &  &  &  &  & 102/3 & \\
 & 32 & -16 &  &  &  & -2 & 1 &  &  &  &  &  &  & 102/3 & \\
64 & -16 &  & 1 &  & 1 &  &  &  &  &  &  &  &  & 102/3 & \\
\hline
& 4 & 1 &  &  &  &  & 1 &  &  &  &  &  &  & 104/1 & \\
\hline
& 12 & -16 &  &  &  &  & 1 &  &  &  &  &  &  & 108/1 & \\
 & 12 &  & -8 &  &  &  & 1 &  &  &  &  &  &  & 108/1 & $I$\\
\hline
 &  &  & 96 &  &  & -9 & 24 &  &  & -16 &  &  &  & 108/2 & $I$\\
\hline
176 & -27 &  &  & 1 &  &  &  &  &  &  &  &  &  & 110/2 & $S$\\
\hline
128 &  & -32 &  &  &  &  & -1 &  &  &  &  &  &  & 110/5 & \\
\hline
 & 64 & -64 & 16 &  &  & -1 &  &  &  &  &  &  &  & 114/4 & $I$\\
\hline
&  &  & 24 &  & 24 & -9 & -40 &  &  &  &  &  &  & 120/1 & \\
 &  & 1 & 1 &  &  &  & -2 & -1 &  & -4 & 5 & -1 &  & 120/1 & \\
 &  & 4 & 4 & 1 &  & -9 & 2 & -1 & -1 &  &  &  &  & 120/1 & \\
 &  & 9 & -1 &  &  &  & -3 &  &  &  &  &  &  & 120/1 & $I$\\
 & 60 & -25 & -36 &  &  & 12 &  &  &  &  &  &  &  & 120/1 & $I$\\
4 & -8 & -1 &  &  &  & 4 & 1 &  &  &  &  &  &  & 120/1 & \\
32 & 44 & -5 &  &  &  & 32 &  &  &  &  &  &  &  & 120/1 & $I$\\
\hline
&  &  & 8 &  & -16 &  &  & 4 &  &  &  & -1 &  & 120/2 & \\
 & 2 &  & -3 &  &  & -2 &  & 5 &  &  &  & -3 &  & 120/2 & \\
 & 3 &  & -2 &  &  & -3 &  & 5 &  &  &  & -2 &  & 120/2 & \\
 & 4 & -4 &  &  &  & -1 & 1 &  &  &  &  &  &  & 120/2 & \\
 & 4 &  & -8 &  &  & 5 & -1 &  &  &  &  &  &  & 120/2 & $I$\\
 & 5 &  & 1 &  & 1 & 4 &  &  &  &  &  &  &  & 120/2 & \\
 & 9 &  & -4 & 1 &  &  &  &  &  &  &  &  &  & 120/2\footnote{[0] Section 4.8.3} & $S$\\
 & 16 &  & -10 & 3 & 2 & -1 &  &  &  &  &  &  &  & 120/2 & \\
 & 16 & 4 &  &  &  & -10 & 6 &  &  & -1 &  &  &  & 120/2 & \\
 & 20 &  & -1 &  & -1 & 1 &  &  &  &  &  &  &  & 120/2 & \\
 & 54 &  & -9 & 1 &  &  &  &  &  &  &  &  &  & 120/2\footnote{[0] Sections 4.8.3, 4.11.6} & $SX_2$\\
 & 81 &  &  &  &  & 27 &  & -16 & 4 &  &  &  &  & 120/2 & \\
 & 216 &  & -36 & 5 &  &  &  &  &  &  &  &  &  & 120/2\footnote{[0] Sections 4.8.3, 4.11.6} & $SX_2$\\
4 & -8 & -1 &  &  &  & 4 &  &  &  &  &  &  &  & 120/2 & \\
4 & -1 & -1 &  &  &  &  &  &  &  &  &  &  &  & 120/2 & $S$\\
16 & -8 & -4 & 1 &  & 1 & 1 &  &  &  &  &  &  &  & 120/2 & \\
16 & -4 & -24 &  &  &  &  & 1 &  &  & 5 &  &  &  & 120/2 & $I$\\
16 & 2 & -4 & 1 &  &  &  &  &  &  &  &  &  &  & 120/2 & $SI$\\
16 & 8 & -8 & -2 & & -2 & & & & & 1 & & & & 120/2 & \\
96 & -12 & 25 & -10 & 1 &  &  &  &  &  &  &  &  &  & 120/2 & $S$\\
\hline
&  &  &  &  & 4 &  & 4 &  &  &  & -15 &  &  & 120/4 & \\
 &  & 1 & -1 &  &  &  &  &  &  &  &  &  &  & 120/4\footnote{[0] Sections 4.4, 4.11.6} & $SI$\\
 &  & 64 &  & -36 &  &  & -16 &  & 9 &  &  &  &  & 120/4 &\\
12 & -24 & 5 &  &  &  & 12 & -5 &  &  &  &  &  &  & 120/ 4 & $I$\\
\hline
 &  & 4 & -6 &  & -6 & 9 &  &  &  &  &  &  &  & 120/5 & \\
 &  & 24 & -2 &  & -6 &  &  & -1 &  &  &  & 1 &  & 120/5 & \\
4 & -8 & -25 &  &  &  & 4 &  &  &  &  &  &  &  & 120/5 & \\
20 & -1 & -9 & 3 & & 3 & -1 & & & & & & & & 120/5 & \\
144 &  & -27 & 15 & -2 &  &  &  &  &  &  &  &  &  & 120/5 & $S$\\
144 & 24 & -16 & 4 & -1 &  &  &  &  &  &  &  &  &  & 120/5 & $SC_6$\\
192 & -24 & 16 & -4 & 1 &  &  &  &  &  &  &  &  &  & 120/5\footnote{[0] Section 4.11.6, [113], [144]} & $SC_6$\\
480 & -12 & 1 &  &  &  &  &  &  &  &  &  &  &  & 120/5 & $SI$\\
\hline
 &  & 8 & -7 &  & -7 & 3 & 5 &  &  & -2 &  &  &  & 126/2 & \\
 &  & 16 & -7 &  & -7 & 3 &  &  &  &  &  &  &  & 126/2 & \\
 & 4 &  & -4 &  & -4 & -1 & 4 &  &  &  &  &  &  & 126/2\footnote{[0] Section 4.8.2} & \\
64 & -28 &  & 12 & -3 &  &  &  &  &  &  &  &  &  & 126/2 & $S$\\
\hline
&  &  & 2 &  & 2 & 1 & -2 &  &  &  &  &  &  & 128/1 & \\
 &  &  & 4 & -1 & -4 & 1 &  &  &  &  &  &  &  & 128/1 & \\
 & 4 & -16 & 3 &  & 3 &  &  &  &  &  &  &  &  & 128/1 & \\
 & 4 &  & -2 &  & -6 &  & 3 &  &  &  &  &  &  & 128/1 & $I$\\
 & 4 &  & -1 &  & -1 &  &  &  &  &  &  &  &  & 128/1 & \\
 & 4 & 2 & -4 & 1 &  &  &  &  &  &  &  &  &  & 128/1\footnote{[0] Section 4.7} & $SC_7$\\
 & 8 & -8 & 2 &  & 2 & -1 &  &  &  &  &  &  &  & 128/1 & \\
8 & -4 & -2 & 4 & -1 &  &  &  &  &  &  &  &  &  & 128/1\footnote{[0] Section 4.7} & $SC_7$\\
16 & -12 & -8 & 2 &  & 2 &  & 5 &  &  & 1 &  &  &  & 128/1 & $I$\\
16 & -8 &  & 2 &  & 2 & -1 &  &  &  &  &  &  &  & 128/1 & \\
16 &  & -12 & 7 & -1 &  &  &  &  &  &  &  &  &  & 128/1 & $S$\\
\hline
& 2 &  & -2 &  & 1 & -1 & -2 & 2 &  &  & 1 & -1 &  & 130/2 & \\
\hline
81 & -9 &  & 1 &  &  &  &  &  &  &  &  &  &  & 132/1 & $SI$\\
\hline
16 & -4 & 1 & 4 & -1 & & -1 & & & & & & & & 140/3 & \\
\hline
 & 48 & 64 &  &  &  & -39 & 12 &  &  &  &  &  &  & 156/2 & $I$\\
\hline
& 4 & 16 &  &  &  &  & -1 &  &  & -4 &  &  &  & 160/1 & \\
 & 16 &  &  &  & -8 & -6 &  & 4 & -1 &  &  &  &  & 160/1 & \\
64 & 4 & -16 &  &  &  &  & -1 &  &  &  &  &  &  & 160/1 & \\
\hline
& 3 & -4 & 1 &  &  &  &  &  &  &  &  &  &  & 162/1 & $SI$\\
\hline
 &  &  &  &  & 36 &  & -60 &  &  &  & 25 & -6 &  & 168/1 & \\
 &  &  & 8 &  & -8 &  & 8 & -4 &  &  & -4 & 3 &  & 168/1 & \\
&  & 108 & -54 &  &  &  &  &  &  & 1 &  &  &  & 168/1 & $I$\\
 & 8 & -4 & -4 &  &  & 1 & -2 &  &  & 1 &  &  &  & 168/1 & $I$\\
 & 8 & -4 &  &  &  & -2 & 1 &  &  &  &  &  &  & 168/1 & \\
 & 8 &  & -1 &  & -1 & -2 & 1 &  &  &  &  &  &  & 168/1 & \\
 & 16 & -16 & 4 & & 4 & -1 & & & & & & & & 168/1 &\\
2 & -2 &  & 1 &  & 1 &  &  &  &  &  &  &  &  & 168/1 & \\
4 & -4 & -9 & 4 &  &  &  &  &  &  &  &  &  &  & 168/1 & $SI$\\
20 & -1 & -5 &  &  &  & -1 & -1 &  &  &  &  &  &  & 168/1 & \\
28 & -4 & 1 &  &  &  &  &  &  &  &  &  &  &  & 168/1 & $S$\\
48 & -6 & 4 & -1 &  &  &  &  &  &  &  &  &  &  & 168/1\footnote{[0] Section 4.11.6} & $SI$\\
96 & -12 & 1 & -2 & 1 &  &  &  &  &  &  &  &  &  & 168/1 & $S$\\
\hline
&  & 4 &  &  &  & 54 & -22 &  &  & -1 &  &  &  & 168/2 & \\
 &  & 8 & -18 &  &  &  & 7 &  &  & -2 &  &  &  & 168/2 & $I$\\
 &  & 16 & -1 &  &  &  &  &  &  & -4 &  &  &  & 168/2 & $I$\\
 & 36 & -4 &  &  &  & -9 & 1 &  &  &  &  &  &  & 168/2 & \\
\hline
1 & 1 &  & -3 &  &  & 1 &  &  &  &  &  &  &  & 180/1 & $I$\\
16 & -4 & -4 & 2 &  & 2 &  & -1 &  &  &  &  &  &  & 180/1 & \\
\hline
 & 4 & -1 & -2 &  &  & 1 &  &  &  &  &  &  &  & 184/1 & $I$\\
92 & -12 & 1 & -2 & 1 &  &  &  &  &  &  &  &  &  & 184/2 & $S$\\
\hline
27 &  & -36 &  &  &  &  & 36 &  &  &  & -32 &  &  & 195/2 & \\
\hline
 & 6 & 16 & -15 & 3 &  &  &  &  &  &  &  &  &  & 198/1 & $S$\\
 & 8 &  & -4 &  &  & 4 &  &  &  &  &  & -1 &  & 198/1 & \\
48 & -9 &  & 1 &  &  &  &  &  &  &  &  &  &  & 198/1 & $SI$\\
48 &  & -24 & -32 &  &  &  &  &  &  & 3 &  &  &  & 198/1 & $I$\\
48 &  & -24 &  &  &  &  & 16 &  &  & 3 &  &  &  & 198/1 & \\
\hline
 & 1 &  & 1 &  &  & -1 &  &  &  &  &  & 1 &  & 210/6 & \\
 & 18 &  & 7 & -3 &  &  &  &  &  &  &  &  &  & 210/6\footnote{[0] Section 4.8.3} & $S$\\
 & 27 &  & -27 & 8 &  &  &  &  &  &  &  &  &  & 210/6\footnote{[0] Section 4.8.3} & $S$\\
\hline
 & 4 & -16 &  &  &  & -1 & 4 &  &  &  &  &  &  & 210/10 & \\
\hline
 &  &  &  &  &  &  & 16 & -6 &  &  & -8 & 3 &  & 216/1 & \\
 &  &  & 1 &  & 1 & -3 & 1 &  &  &  &  &  &  & 216/1 & \\
 &  &  & 4 &  & 4 & -3 &  &  &  &  &  &  &  & 216/1 & \\
 & 12 &  & -1 &  & -1 &  &  &  &  &  &  &  &  & 216/1 & \\
 & 12 &  & 2 &  & 2 &  & -1 &  &  &  &  &  &  & 216/1 & \\
 & 16 &  & -12 & 3 & 4 & -1 &  &  &  &  &  &  &  & 216/1 & \\
 & 32 &  &  & -3 & -16 & 7 &  &  &  &  &  &  &  & 216/1 & \\
 & 64 &  &  & -4 &  & -16 &  &  & 1 &  &  &  &  & 216/1 & \\
 & 64 &  &  &  &  & -28 &  & 12 & -3 &  &  &  &  & 216/1 & \\
 & 256 &  &  &  &  & 32 &  &  & -3 &  &  &  &  & 216/1 & \\
\hline
 &  &  &  &  & 4 &  & -12 & 2 &  &  & 9 & -3 &  & 216/2 & \\
 &  &  &  &  & 16 &  & 16 &  &  &  & 4 &  & -1 & 216/2 & \\
 & 64 &  &  & -3 &  & -19 &  & 3 &  &  &  &  &  & 216/2 & \\
\hline
 & 27 & -4 & -27 &  &  &  & 4 &  &  &  &  &  &  & 258/1 & $I$\\
\hline
 &  & 27 & -27 &  &  &  & 27 &  &  & -4 &  &  &  & 264/2 & $I$\\
 & 32 &  & 8 &  & 8 & -11 &  &  &  &  &  &  &  & 264/2 & \\
\hline
 &  & 24 & -2 &  &  &  & -9 &  &  & -6 &  &  &  & 264/4 & $I$\\
 & 1 &  & -2 &  &  & -1 &  & 3 &  &  &  & -2 &  & 264/4 & \\
 & 6 &  & -1 &  &  & -6 &  & 7 &  &  &  & -1 &  & 264/4 & \\
 & 18 &  & -11 & 3 &  &  &  &  &  &  &  &  &  & 264/4\footnote{[0] Section 4.8.3} & $S$\\
 & 27 &  &  & -1 &  &  &  &  &  &  &  &  &  & 264/4\footnote{[0] Section 4.8.3} & $S$\\
\hline
 &  &  & 4 &  & 4 & 1 &  &  &  &  &  &  &  & 280/2 & \\
 & 16 &  & -16 & 5 &  &  &  &  &  &  &  &  &  & 280/2\footnote{[0] Section 4.8.3} & $S$\\
 & 32 &  &  & 1 & -16 & 3 &  &  &  &  &  &  &  & 280/2 & \\
 & 64 &  &  &  &  & -12 &  & -4 & 1 &  &  &  &  & 280/2 & \\
\hline
 &  &  &  &  &  &  & 24 & 16 &  &  & -12 & -5 & 2 & 288/1 & \\
 &  &  & 1 &  & -1 &  &  &  &  &  &  &  &  & 288/1 & $I$\\
 &  &  & 2 &  & 2 & 3 & -2 &  &  &  &  &  &  & 288/1 & \\
 & 4 & -6 & 2 &  & 2 & -1 &  &  &  &  &  &  &  & 288/1 & \\
 & 4 &  & -4 &  & -4 & -1 & 3 &  &  &  &  &  &  & 288/1 & \\
 & 4 &  & -1 &  & -1 & -1 &  &  &  &  &  &  &  & 288/1 & \\
 & 12 &  & -4 &  & -4 &  & 1 &  &  &  &  &  &  & 288/1 & \\
8 & -4 & -2 & 2 &  & 2 & -1 &  &  &  &  &  &  &  & 288/1 & \\
16 & -16 & -4 & 12 & -3 &  &  &  &  &  &  &  &  &  & 288/1\footnote{[0] Section 4.11.7} & $S$\\
16 & -4 & -8 & 2 & & 2 & & 3 & & & 1 & & & & 288/1 & \\
16 & -4 & 4 & -1 &  & -1 &  &  &  &  &  &  &  &  & 288/1 & \\
48 &  & -12 & 7 & -1 &  &  &  &  &  &  &  &  &  & 288/1 & $S$\\
\hline
3 & 6 & -12 & 18 &  &  &  & -20 &  &  &  &  &  &  & 300/2 & $I$\\
\hline
 & 16 &  &  &  & 16 & 8 &  & -8 & 3 &  &  &  &  & 330/3 & \\
\hline
 & 256 &  & -128 &  &  & -49 &  &  &  &  &  &  &  & 336/1 & $I$\\
\hline
 & 2 &  & -3 & 1 &  & -2 &  & 5 & -1 &  &  & -3 & 1 & 360/2 & \\
 & 6 &  & -1 &  &  &  &  &  &  &  &  &  &  & 360/2\footnote{[0] Sections 4.8.3, 4.11.6} & $SIX_1$\\
 & 12 & -16 &  &  &  &  & -3 &  &  & 4 &  &  &  & 360/2 & \\
 & 12 &  & 16 &  &  &  & -3 &  &  &  &  &  &  & 360/2 & $I$\\
 & 24 &  & -4 & 1 &  &  &  &  &  &  &  &  &  & 360/2\footnote{[0] Sections 4.8.3, 4.11.6} & $SX_1$\\
 & 48 & -24 &  &  &  &  & -2 &  &  & 1 &  &  &  & 360/2 & \\
96 &  &  & -1 &  &  &  &  &  &  &  &  &  &  & 360/2 & $SI$\\
\hline
 &  & 4 & -4 &  &  & -1 & 2 &  &  & -1 &  &  &  & 380/1 & $I$\\
\hline
 &  & 216 & -27 &  &  &  &  &  &  & -4 &  &  &  & 384/1 & $I$\\
 & 4 &  & -2 &  & -2 & 2 & 3 & -3 &  &  & -1 & 1 &  & 384/1 & \\
 & 4 &  & 4 &  &  & -4 & 1 &  &  &  &  &  &  & 384/1 & $I$\\
8 & -8 & -2 &  &  &  & 3 & -1 &  &  &  &  &  &  & 384/1 & \\
16 & 8 &  & 2 &  & 2 & 3 &  &  &  &  &  &  &  & 384/1 & \\
64 & -4 & -16 &  &  &  &  & 1 &  &  &  &  &  &  & 384/1 & \\
64 & -4 &  & 1 &  & 1 &  &  &  &  &  &  &  &  & 384/1 & \\
128 & -20 &  & 1 &  & 1 &  &  &  &  &  &  &  &  & 384/1 & $I$\\
\hline
 & 108 & -16 &  &  &  &  & -11 &  &  & 4 &  &  &  & 384/3 & $I$\\
\hline
 &  & 8 & 23 & -3 & -2 & -6 &  &  &  &  &  &  &  & 390/5 & \\
 & 8 &  & -108 & 27 &  & 25 &  &  &  &  &  &  &  & 390/5 & \\
\hline
 & 4 & -4 & -2 &  & -2 &  & -3 &  &  & 1 &  &  &  & 392/2 & \\
\hline
 &  &  &  &  & 24 &  & -24 &  &  &  & 6 & -1 &  & 396/1 & \\
 &  & 2 & 1 &  & -1 &  & 1 &  &  &  &  &  &  & 396/1 & $I$\\
 & 8 &  & 1 &  & 1 & -2 & -1 &  &  &  &  &  &  & 396/1 & \\
 & 24 & 4 & -8 &  & -8 &  & 4 &  &  & -1 &  &  &  & 396/1 & $I$\\
4 & 2 &  & 1 &  &  &  &  &  &  &  &  &  &  & 396/1 & $SI$\\
48 &  & -12 &  &  &  &  & 1 &  &  &  &  &  &  & 396/1 & \\
\hline
128 &  & -32 &  &  &  &  & 6 &  &  &  & -1 &  &  & 416/2 & \\
\hline
44 & -3 & 25 & -20 & 4 &  &  &  &  &  &  &  &  &  & 440/2 & $S$\\
\hline
 &  & 4 & -81 &  &  &  & 32 &  &  &  &  &  &  & 480/1 & $I$\\
\hline
 &  &  &  &  & 4 &  & 16 &  &  &  & 15 &  &  & 480/2 & \\
 &  & 4 & -7 & 3 &  &  &  &  &  &  &  &  &  & 480/2\footnote{[0] Sections 4.4, 4.11.6} & $S$\\
 &  & 108 & 27 &  &  &  &  &  &  & -32 &  &  &  & 480/2 & $I$\\
 & 4 &  & 2 &  & -2 & 2 & 1 & -1 &  &  & 1 & -1 &  & 480/2 & \\
 & 8 &  & -1 &  & -1 & 4 & -1 &  &  &  &  &  &  & 480/2 & \\
 & 8 &  & 2 &  & 2 & -5 & -4 &  &  &  &  &  &  & 480/2 & \\
\hline
 &  &  &  &  & 12 &  & 16 &  &  &  & -3 &  &  & 480/5 & \\
 &  &  &  &  & 27 &  & -54 &  &  & -50 & 27 &  &  & 480/5 & \\
 &  & 12 & -13 & 3 &  &  &  &  &  &  &  &  &  & 480/5\footnote{[0] Sections 4.4, 4.11.6} & $S$\\
 & 36 &  & -6 &  & -30 &  & 5 &  &  &  &  &  &  & 480/5 & $I$\\
8 & -16 & -2 &  &  &  & 8 & -3 &  &  &  &  &  &  & 480/5 & \\
72 &  & -18 &  &  &  &  & 1 &  &  &  &  &  &  & 480/5 & \\
72 &  & -10 & 3 &  &  &  &  &  &  &  &  &  &  & 480/5 & $SI$\\
\hline
48 & -12 & -3 &  &  &  &  & 1 &  &  &  &  &  &  & 516/1 & \\
\hline
 &  & 20 & -20 &  &  & 3 & -6 &  &  & 3 &  &  &  & 540/2 & $I$\\
\hline
 &  & 48 &  &  &  & -1 & -4 &  &  & -12 &  &  &  & 544/1 & \\
 & 4 & -12 & 7 & -1 & 3 & -3 & & & & & & & & 544/1 & \\
\hline
 & 8 & -4 &  &  &  & -2 & -1 &  &  & 1 &  &  &  & 552/1 & \\
 & 8 & 4 & -2 &  & -2 & 1 & 2 &  &  &  &  &  &  & 552/1 & \\
16 & 8 & -4 & 1 &  & 1 &  &  &  &  &  &  &  &  & 552/1 & \\
\hline
 & 12 & -1 & 2 & -1 &  &  &  &  &  &  &  &  &  & 552/2 & $S$\\
 & 64 &  &  & -27 &  & -43 &  & 27 &  &  &  &  &  & 552/2 & \\
\hline
16 & 13 &  & -45 &  &  & 16 &  &  &  &  &  &  &  & 570/7 & $I$\\
\hline
 &  & 16 &  &  &  & 54 & -25 &  &  &  &  &  &  & 606/1 & \\
\hline
 & 54 & 16 & -27 &  &  &  & 2 &  &  &  &  &  &  & 618/1 & $I$\\
\hline
 & 24 & 4 & -4 &  &  & -9 & 6 &  &  & -1 &  &  &  & 648/1 & $I$\\
\hline
 & 12 & 16 &  &  &  & -3 &  &  &  &  &  &  &  & 678/2 & \\
 & 18 & -16 & 3 &  &  &  &  &  &  &  &  &  &  & 678/2 & $SI$\\
\hline
 &  & 64 &  &  &  & -3 & -24 &  &  &  &  &  &  & 684/4 & \\
\hline
 & 1 & -27 & 3 &  & 3 &  & 5 &  &  &  &  &  &  & 760/1 & \\
\hline
 &  &  &  &  &  &  & 8 & -24 & 18 &  & -4 & -3 &  & 864/1 & \\
\hline
 & 36 &  & -2 &  & -18 & 3 & 1 &  &  &  &  &  &  & 918/1 & $I$\\
\hline
16 &  & -4 & 1 &  &  &  &  &  &  &  &  &  &  & 928/1 & $SI$\\
\hline
 & 16 & -2 & 3 &  & -27 &  & -5 &  &  &  &  &  &  & 940/1 & $I$\\
\hline
16 & -2 & -4 & 1 &  &  &  &  &  &  &  &  &  &  & 1160/1 & $SI$\\
\hline
 &  & 16 &  &  &  & -3 & -12 &  &  & -4 &  &  &  & 1440/7 & \\
\hline
 & 12 & 4 &  &  &  & -3 & 3 &  &  &  &  &  &  & 1608/1 & \\
4 & -8 &  & 3 &  &  & 4 &  & -3 &  &  & 1 &  &  & 1608/1 & \\
\hline
 & 18 & -16 & 4 &  & -1 &  &  &  &  &  &  &  &  & 1890/4 & $I$\\
\hline
 & 36 & -32 & 3 &  & 3 &  &  &  &  &  &  &  &  & 1920/3 & \\
 & 36 & -16 &  &  &  &  & -9 &  &  & 4 &  &  &  & 1920/3 & \\
 & 36 &  & -12 &  &  &  & 1 &  &  &  &  &  &  & 1920/3 & $I$\\
\end{longtable}
\end{scriptsize}
\end{center}

\section{References}

This paper can be read as apocrypha to my book:
\begin{itemize}
\item[{[0]}] Meyer, Christian, \textit{Modular Calabi-Yau threefolds},
Fields Institute Monograph \textbf{22} (2005), AMS.
\end{itemize}

\subsection{Updated references}

The following papers are contained as preprints in the references of [0]. Meanwhile they have been published or accepted for publication. In addition, the URL of the Modular Forms Database ([97]) has changed (repeatedly).

\begin{itemize}
\item[{[13]}] Bernardara,~M.,
\textit{Calabi--Yau complete intersections with infinitely many lines},
Rend. Sem. Mat. Univ. Pol. Torino \textbf{66}, no. 2 (2008), pp. 87--97.

\item[{[28]}] Cynk,~S., Meyer,~C.,
\textit{Geometry and arithmetic of certain double octic Calabi--Yau manifolds},
Canadian Math. Bull. \textbf{48}, no. 2 (2005), pp. 180--194.

\item[{[30]}] Cynk,~S., van~Straten,~D.,
\textit{Infinitesimal deformations of double covers of smooth algebraic varieties},
Math. Nachr. \textbf{279}, no. 7 (2006), pp. 716--726.

\item[{[34]}] Dieulefait,~L.,
\textit{A modularity criterion for integral Galois representations and Calabi--Yau threefolds},
AMS/IP Studies in Advanced Mathematics ''Mirror Symmetry V'' (2006),
Proceedings of the BIRS Workshop on Calabi-Yau Varieties and Mirror Symmetry, December 6-11, 2003,
appendix to Hulek,~K., Verrill,~H.~A.,
\textit{On the modularity of Calabi--Yau threefolds containing elliptic ruled surfaces}, pp. 32--34.

\item[{[51]}] Hulek,~K., Verrill,~H.~A.,
\textit{On modularity of rigid and nonrigid Calabi--Yau varieties associated to the root lattice $A_4$},
Nagoya Math. Journal \textbf{179} (2005), pp. 103--146.

\item[{[52]}] Hulek,~K., Verrill,~H.~A.,
\textit{On the modularity of Calabi--Yau threefolds containing elliptic ruled surfaces},
AMS/IP Studies in Advanced Mathematics ''Mirror Symmetry V'' (2006),
Proceedings of the BIRS Workshop on Calabi-Yau Varieties and Mirror Symmetry, December 6-11, 2003, pp. 19--31.

\item[{[55]}] Kimura,~K.,
\textit{A rational map between two threefolds},
AMS/IP Studies in Advanced Mathematics ''Mirror Symmetry V'' (2006),
Proceedings of the BIRS Workshop on Calabi-Yau Varieties and Mirror Symmetry, December 6-11, 2003, pp. 87--88.

\item[{[64]}] Livn\'e,~R., Yui,~N.,
\textit{The modularity of certain non-rigid Calabi--Yau threefolds},
J. of Math. Kyoto Univ. \textbf{45}, no. 4 (2005), pp. 645--665.

\item[{[74]}] Mortenson,~E., \textit{Modularity of a Certain Calabi--Yau Threefold and
Combinatorial Congruences}, Ramanujan Journal \textbf{11}, no. 1 (2006), pp. 5--39.

\item[{[89]}] Sch\"utt,~M.,
\textit{On the modularity of three Calabi--Yau threefolds with bad reduction at 11},
Canadian Math. Bull. \textbf{49}, no. 2 (2006), pp. 296--312.

\item[{[97]}] Stein,~W.~A.,
\textit{The Modular Forms Database} (2005),\\
{\tt https://wstein.org/Tables}.

\item[{[98]}] Top,~J., van Geemen,~B.,
\textit{An isogeny of K3 surfaces},
Bull. London Math. Soc. \textbf{38} (2006), pp. 209--223.

\item[{[105]}] Wan,~D.,
\textit{Mirror Symmetry for Zeta Functions}, with an appendix by C. D. Haessig,
AMS/IP Studies in Advanced Mathematics ''Mirror Symmetry V'' (2006),
Proceedings of the BIRS Workshop on Calabi-Yau Varieties and Mirror Symmetry, December 6-11, 2003, pp. 159--184.
\end{itemize}

\subsection{Additional references}

Most of the following papers have appeared after [0] was in print. Most of them are dealing with modularity of Calabi-Yau varieties in some way. References from [152] are new in this version of the paper.

The paper [146] is a classic. Only after publication of [0], I discovered that in this paper the family of $\Sigma_5$-symmetric quintic surfaces (cf. [0], Section 4.9) has been studied in detail. It is proved that the quintic surface $S_{(5:12)}$ is a singular model for the Hilbert modular surface corresponding to the congruence subgroup of level 2 in the extended Hilbert modular group for $\qz(\sqrt{13})$. This provides an explanation for the occurrence of the bad prime 13 in the $L$-series of the double octic $X_{(5:12)}$ constructed from $S_{(5:12)}$ and the Clebsch cubic. Furthermore, the surface $S_{(5:4)}$ corresponds to $\qz(\sqrt{21})$ in the
same way. I do not know, however, how the bad prime 7 gets lost in the level of the $L$-series of the double octic $X_{(5:4)}$.

For theoretical advances on modularity of Calabi-Yau threefolds, cf. [129]. For \emph{really} interesting examples I recommend [114] and [138].

\begin{itemize}

\item[{[113]}] Bini,~G., van Geemen,~B., \textit{Geometry and Arithmetic of Maschke's Calabi-Yau Threefold},
Comm. in Number Theory and Physics \textbf{5}, no. 4 (2011), pp. 779--826.

\item[{[114]}] Brown,~F., Schnetz,~O., \textit{Modular forms in quantum field theory},
Comm. in Number Theory and Physics \textbf{7}, no. 2 (2013), pp. 293--325.

\item[{[115]}] Burek,~D., \textit{Rigid realizations of modular forms in Calabi--Yau threefolds},
J. Pure and Applied Algebra \textbf{223}, no. 2 (2019), pp. 547--552.

\item[{[116]}] Cynk,~S., Kocel-Cynk,~B., \textit{Classification of double octic Calabi-Yau threefolds},
preprint (2016), arXiv: math.AG/1612.04364.

\item[{[117]}] Cynk,~S., Freitag,~E., Salvato Manni,~R., \textit{The geometry and arithmetic of a Calabi-Yau Siegel threefold},
International J. Math. \textbf{22}, no. 11 (2011), pp. 1585--1602.

\item[{[118]}] Cynk,~S., Hulek,~K., \textit{Construction and examples of higher-dimen\-sional modular Calabi-Yau manifolds},
Canadian Math. Bull. \textbf{50}, no. 4 (2007), pp. 486--503.

\item[{[119]}] Cynk,~S., Meyer,~C., \textit{Modular Calabi-Yau threefolds of level eight},
International J. Math. \textbf{18}, no. 3 (2007), pp. 331--347.

\item[{[120]}] Cynk,~S., Meyer,~C., \textit{Modularity of some non-rigid double octic Cala\-bi-Yau threefolds},
Rocky Mountain J. Math. \textbf{38}, no. 6 (2008), pp. 1937--1958.

\item[{[121]}] Cynk, S., Sch\"utt, M., \textit{Generalised Kummer constructions and Weil restrictions},
J. Number Theory \textbf{129}, no. 8 (2009), pp. 1965--1975 .

\item[{[122]}] Cynk, S., Sch\"utt, M., \textit{A modular Calabi-Yau threefold with CM by $\qz(\sqrt{-23})$},
supplement to Cynk, S., Sch\"utt, M., \textit{Generalised Kummer constructions and Weil restrictions} (2009),\\
{\tt http://www.math.ku.dk/$\sim$mschuett/publik\_en.php.html}

\item[{[123]}] Dieulefait,~L., \textit{Improvements on Dieulefait-Manoharmayum and applications},
preprint (2005), arXiv: math.NT/0508163.

\item[{[124]}] Dieulefait, L., \textit{On the modularity of rigid Calabi-Yau threefolds: Epilogue},
Proceedings of the trimester on Diophantine Equations at the Hausdorff Institute, Zapiski POMI of the Steklov Math. Inst. St. Petersbourg \textbf{377} (2010), pp. 44--49,
also in J. Math. Sciences. \textbf{171} (2010), pp. 725--727.

\item[{[125]}] Dieulefait,~L., Pacetti,~A., Sch\"utt,~M., \textit{Modularity of the Consani-Scholten quintic. With an appendix by Jos\'e Burgos Gil and Ariel Pacetti},
Doc. Math. \textbf{17} (2012), pp. 953--987.

\item[{[126]}] Elkies, N., Sch\"utt, M., \textit{Modular forms and K3 surfaces}, Adv. Math. \textbf{240}, no. 20 (2013), pp. 106--131.

\item[{[127]}] Fu,~L., Wan,~D., \textit{Mirror congruence for rational points on Calabi-Yau varieties},
Asian J. Math. \textbf{10} (2006), pp. 1--10.

\item[{[128]}] Gouv\^{e}a,~F., Kiming,~I., Yui,~N., \textit{Quadratic twists of rigid Calabi-Yau threefolds over $\qz$},
in \textit{Arithmetic and Geometry of K3 Surfaces and Calabi–Yau Threefolds}, eds. Laza,~R., Schütt,~M., Yui,~N. (2013), Fields Institute Communications \textbf{67}, Springer, pp. 517--533.

\item[{[129]}] Gouv\^{e}a, F., Yui, N., \textit{Rigid Calabi-Yau Threefolds over $\qz$ are Modular},
Expo. Math. \textbf{29}, no. 1 (2011), pp. 142--149.

\item[{[130]}] Hulek,~K., Kloosterman,~R., \textit{The L-series of a cubic fourfold},
Manu\-scripta Math. \textbf{124}, no. 3 (2007), pp. 391--407.

\item[{[131]}] Hulek,~K., Kloosterman,~R., Sch\"utt,~M., \textit{Modularity of Calabi-Yau varieties},
in \textit{Global Aspects of Complex Geometry}, eds. Catanese,~F., Esnault,~H., Huckleberry,~A., Hulek,~K., Peternell,~T. (2006),
Sprin\-ger, pp. 271--309.

\item[{[132]}] Hulek,~K., Verrill,~H.~A., \textit{On the motive of Kummer varieties associated to $\Gamma_1(7)$ - Supplement to the
paper: The modularity of certain non-rigid Calabi-Yau threefolds (by R. Livn\'e and N. Yui)},
J. of Math. Kyoto Univ. \textbf{45}, no. 4 (2005), pp. 667--681.

\item[{[133]}] Kapustka,~G., Kapustka,~M., \textit{Modularity of a nonrigid Calabi--Yau manifold with bad reduction at 13},
Ann. Polon. Math. \textbf{90} (2007), pp. 89--98.

\item[{[134]}] Kapustka,~M., \textit{Correspondences between modular Calabi--Yau fiber products},
Manuscripta Math. \textbf{130}, no. 1 (2009), pp. 121--135.

\item[{[135]}] Lee,~E., \textit{A Modular Non-Rigid Calabi-Yau Threefold},
AMS/IP Studies in Advanced Mathematics ''Mirror Symmetry V'' (2006),
Proceedings of the BIRS Workshop on Calabi-Yau Varieties
and Mirror Symmetry, December 6-11, 2003, pp. 89--122.

\item[{[136]}] Lee,~E., \textit{A modular quintic Calabi-Yau threefold of level 55}
Canad. J. Math. \textbf{63} (2011), pp. 616--633.

\item[{[137]}] Lee,~E., \textit{Update on Modular Non-Rigid Calabi-Yau Threefolds},
in \textit{Modular Forms and String Duality} (2008),
eds. Yui, N., Verrill, H., Doran, C.F., Fields Inst. Commun. Series, AMS, pp. 65--82.

\item[{[138]}] Logan,~A., \textit{New realizations of modular forms in Calabi–Yau three\-folds arising from $\phi^4$ theory},
J. Number Theory \textbf{184}, no. 3 (2018), pp. 342--383 .

\item[{[139]}] Meyer,~C., \textit{The mirror quintic as a quintic},
preprint (2005), arXiv: math.AG/0503329.

\item[{[140]}] Samol, K., van Straten, D., \textit{Frobenius polynomials for Calabi-Yau equations},
Comm. in Number Theory and Physics \textbf{2}, no. 3 (2008), pp. 537--561.

\item[{[141]}] Sch\"utt,~M.,
\textit{Arithmetic of a singular K3 surface},
Michigan Math. J. \textbf{56}, no. 3 (2008), pp. 513--527.

\item[{[142]}] Sch\"utt, M.,
\textit{Arithmetic of K3 surfaces},
Jahresbericht der DMV \textbf{111}, no. 1 (2009), pp. 23--41.

\item[{[143]}] Sch\"utt,~M.,
\textit{CM newforms with rational coefficients},
Ramanujan J. \textbf{19} (2009), pp. 187--205.

\item[{[144]}] Sch\"utt, M.,
\textit{Modularity of Maschke's octic and Calabi-Yau threefold},
Comm. in Number Theory and Physics \textbf{5}, no. 4 (2011), pp. 827--847.

\item[{[145]}] Sch\"utt,~M., Top,~J.,
\textit{Arithmetic of the [19,1,1,1,1,1] fibration},
Commentarii Mathematici Universitatis Sancti Pauli \textbf{55},
no. 1 (2006), pp. 9--16.

\item[{[146]}] Van der Geer,~G., Zagier,~D.,
\textit{The Hilbert Modular Group for the Field $\qz(\sqrt{13})$},
Invent. Math. \textbf{42} (1977), pp. 93--133. 

\item[{[147]}] Ward,~M., \textit{Arithmetic Properties of the Derived Category for Calabi-Yau Varieties}, thesis (2014), Univ. of Washington.

\item[{[148]}] Yang, L., \textit{Galois representations arising from twenty-seven lines on a cubic surface and the arithmetic associated with Hessian polyhedra},
preprint (2006), arXiv: math.NT/0612383.

\item [{[149]}] Yui,~N., \textit{Modularity of Calabi-Yau Varieties: 2011 and Beyond},
in \textit{Arithmetic and Geometry of K3 Surfaces and Calabi–Yau Threefolds}, eds. Laza,~R., Schütt,~M., Yui,~N. (2013), Fields Institute Communications \textbf{67}, Springer, pp. 101--139.

\item[{[150]}] Wan,~D., \textit{Arithmetic mirror symmetry},
Quarterly Journal of Pure and Applied Math. \textbf{1}, no. 2 (2005), pp. 369--378.

\item[{[151]}] Zudilin,~W., \textit{A Hypergeometric Version of the Modularity of Rigid Calabi-Yau Manifolds},
SIGMA \textbf{14} (2018), 086.

\item[{[152]}] Bastian,~B., Van de Heisteeg,~D., Schlechter,~L., \textit{Beyond Large Complex Structure: Quantized Periods and Boundary Data for One-Modulus Singularities}, preprint (2023), arXiv:2306.01059.

\item[{[153]}] Bönisch,~K, Klemm,~A., Scheidegger,~E., Zagier,~D., \textit{D-brane masses at special fibres of hypergeometric families of Calabi-Yau threefolds, modular forms, and periods}, preprint (2022), arXiv:2202.09426.

\item[{[154]}] Candelas,~P., De la Ossa,~X., McGovern,~J., \textit{Classical Weight-Four L-value Ratios as Sums of Calabi–Yau Invariants}, preprint (2025), arXiv:2410.07107.

\item[{[155]}] Chmiel,~T., Cynk,~S., \textit{Periods of Singular Double Octic Calabi-Yau Threefolds and Modular Forms}, Math. Nachr. 296 (2023), no. 8, 3257-3271.

\item[{[156]}] Cynk,~S., Kocel-Cynk,~B., \textit{Classification of Double Octic Calabi-Yau Threefolds Defined by an Arrangement of Eight Planes II}, preprint (2026), arXiv:2602.19413.

\item[{[157]}] Donlagic,~A., textit{Numerical Calculation of Periods on Schoen's Class of Calabi-Yau Threefolds}, preprint (2025), arXiv:2504.09383.

\item[{[158]}] Dummigan,~N., \textit{Modularity of a Certain "Rank-2 Attractor" Calabi-Yau Threefold}, preprint (2026), arXiv:2602.20188.

\item[{[159]}] Festi,~D., Van Geemen,~B., \textit{A Calabi-Yau Threefold Coming from Two Black Holes}, preprint (2022), arXiv:2207.01936.

\item[{[160]}] Grzelakowski, K. \textit{On Quintic Threefolds With Triple Point},\\ preprint (2020), arXiv:2010.13866.

\item[{[161]}] Kachru,~S., Nally,~R., Yang,~W, \textit{Supersymmetric Flux Compactifications and Calabi-Yau Modularity}, preprint (2020), arXiv:2001.06022.

\item[{[162]}] Kachru,~S., Nally,~R., Yang,~W, \textit{Flux Modularity, F-Theory, and Rational Models}, preprint (2020), arXiv:2010.07285.

\end{itemize}

\end{document}